\documentclass[12pt]{article}
\pdfoutput=1
\usepackage{amsmath}
\usepackage{url}
\usepackage{graphicx}
\usepackage{epsfig,subfigure}
\usepackage{fullpage}

\newcommand{\bfd}{\mathbf{d}}
\newcommand{\bfdo}{\mathbf{d}_{\mathrm{obs}}}
\newcommand{\bfdto}{\tilde{\mathbf{d}}_{\mathrm{obs}}}
\newcommand{\bfdp}{\mathbf{d}_{\mathrm{pre}}}
\newcommand{\bfe}{{\mathbf e}}

\newcommand{\bfm}{\mathbf{m}}

\newcommand{\bfma}{\mathbf{m}_{\mathrm{apr}}}
\newcommand{\bfme}{\mathbf{m}_{\mathrm{exact}}}

\newcommand{\bfp}{\mathbf{p}}

\newcommand{\bfr}{\mathbf{r}}
\newcommand{\bfrte}{\tilde{\mathbf{r}}_{\mathrm{exact}}}

\newcommand{\bfv}{\mathbf{v}}

\newcommand{\bfu}{\mathbf{u}}
\newcommand{\bfy}{\mathbf{y}}
\newcommand{\bfz}{\mathbf{z}}

\newcommand{\bfze}{\mathbf{z}_{\mathrm{exact}}}

\newcommand{\Wd}{W_{\mathbf d}}

\newcommand{\alphao}{\alpha_{\mathrm{opt}}}

\newcommand{\bfmalpha}{\bfy(\alpha)}

\newcommand{\bfmyalpha}{\bfy}

\newcommand{\argmin}[1]{\textnormal{arg} \min_{#1}}

\newcommand{\figdirA}

\begin{document}
\title{Application of the $\chi^2$ principle and unbiased predictive risk estimator  for determining the regularization parameter in 3D focusing gravity inversion}
\author{Saeed Vatankhah$^1$, Vahid E Ardestani$^1$ and Rosemary A Renaut$^2$\\
$^1$Institute of Geophysics, University of  Tehran, Tehran, Iran\\
$^2$ School of Mathematical and Statistical Sciences, Arizona State University, Tempe, USA.}
\date{\today}
\maketitle

\begin{abstract}
The $\chi^2$ principle and the  unbiased predictive risk estimator are used to determine  optimal regularization parameters in the context of 3D focusing gravity inversion with the minimum support stabilizer. At each iteration of the focusing inversion  the  minimum support stabilizer is  determined and  then the fidelity term is updated using the standard form transformation. Solution of the resulting Tikhonov functional is  found efficiently using the singular value decomposition of the transformed model matrix, which also provides for efficient determination of the updated regularization parameter each step.  Experimental 3D simulations using  synthetic data of a  dipping dike and a cube anomaly  demonstrate that both parameter estimation techniques  outperform the Morozov discrepancy principle for determining the regularization parameter. Smaller  relative errors of the   reconstructed models are obtained with fewer iterations.  Data acquired over the Gotvand dam site in the south-west of Iran are used to validate use of  the methods for inversion of  practical data and provide good estimates of anomalous structures within the subsurface. 
\end{abstract}


\section{Introduction}
Gravity surveys have been used for many years for a  wide range of studies including oil and gas exploration, mining applications, mapping bedrock topography, estimation of the crustal thickness and recently-developed microgravity investigations \cite{nabighian:2005}. The inversion of gravity data is one of the important steps in the interpretation of practical data. The goal is to estimate density and geometry parameters of an unknown subsurface model from a set of known gravity observations measured on the surface. In the linear inversion of gravity data it is standard to assume that the subsurface under the survey area can be approximated through a discretization of the subsurface into rectangular blocks of constant density \cite{Bou:2001}. In solving for the densities at these blocks this kind of parameterization is flexible 
for the  reconstruction of the subsurface model, but requires  more unknowns than observations and thus introduces  algebraic ambiguity in the solution of the linear system. Additionally,  the existence of noise in the measurements of practical data and the inherent non-uniqueness of the gravity sources, based on Gauss's theorem, means that the inversion of gravity data is an example of an underdetermined and  ill-posed problem.  Thus, in order to find an acceptable solution  which is less sensitive to the measurement error regularization, also known as  stabilization, is typically imposed. A popular approach  uses the minimization of a cost functional that combines the data fidelity with an L$2$, or Tikhonov, type regularization, see e.g. \cite{ABT,Regtools,Vogel:2002}. Two important aspects of the Tikhonov regularization are the choices of the stabilizing operator and the regularization parameter. The former  impacts the class of solution which will be obtained, and the latter controls the trade off between the data fit and the regularization term. Two main classes of stabilizer have been used in the inversion of gravity data;  a smoothing stabilizer which  employs the first or second derivative of the model parameters see e.g. \cite{Li:96,Bou:2001} and a stabilizer which produces non-smooth models e.g. \cite{Bou:2001,Last:83,PoZh:99}. In this paper the minimum support (MS) stabilizer which was introduced in \cite{Last:83} and developed in \cite{PoZh:99} is used to reconstruct models with non-smooth features. 

The determination of an optimal regularization parameter in potential field data inversion is a topic of much previous research and includes methods such as the L-curve (LC) \cite{Li:99,Far:2004,vatan:2014}, generalized cross validation (GCV) \cite{Far:2004,vatan:2014} and the more often adopted Morozov discrepancy principle  (MDP) \cite{morozov,Li:96,Far:2004}. 
Because it is well-know that the MDP generally  overestimates the regularization parameter,  hence leading to overly smoothed solutions, we discuss here regularization parameter estimation in the specific context of the inversion of underdetermined gravity data using  the   Unbiased Predictive Risk Estimator (UPRE) and the $\chi^2$ principle,  see e.g.  \cite{Vogel:2002,VRA2014}.
 Whereas in \cite{vatan:2014} we considered the use of the GCV and LC methods for 2D focusing inversion,  our subsequent investigations in \cite{VRA2014} demonstrated that for small scale 2D problems the UPRE  and $\chi^2$ principle   improve on results using the  LC, GCV and MDP, with respect to reduced relative error, reduced computational cost or both. Indeed,  all  methods demonstrated their efficiency  as compared with the MDP  \cite{VRA2014}, but the UPRE and $\chi^2$ techniques offer the most promise for parameter estimation in terms of cost and accuracy. We, therefore,  solve the underlying regularized model, with these parameter-choice methods, here contrasting for completeness with the MDP. Moreover,   in place of the  use of the generalized singular value decomposition (GSVD), \cite{PaigeSau1},   as advocated in \cite{vatan:2014,VRA2014},  we use  the singular value decomposition (SVD) of the system matrix in standard form \cite{Hansen}. This provides a more efficient tool as compared to the GSVD for the solution of Tikhonov regularized problems of small to moderate scale. 

The outline of this paper is as follows. In section~\ref{3D modelling} we review the derivation of the analytic calculation of the gravity anomaly derived from a 3D cell model. In section~\ref{focusing} the algorithm for focusing inversion is discussed.  Furthermore,  in this section numerical solutions of the Tikhonov objective function using the SVD for the regularized-modified model system are discussed. Extensions of the MDP, UPRE and $\chi^2$ methods  for estimating the regularization parameter  have been extensively  discussed in \cite{VRA2014}, but we provide a brief rationale for the latter two methods which are not well-known in this field in section~\ref{parameter estimation} with necessary formulae collected in  \ref{regparam}.   Results for synthetic examples are illustrated in section~\ref{synthetic}. The approach is applied on gravity data acquired from Gotvand dam site in section~\ref{real}. Conclusions and a discussion of future plans follow in section~\ref{conclusion}.

\section{Gravity modelling}\label{3D modelling}
\begin{figure}
\label{3DModel}
\includegraphics[width=.9\textwidth]{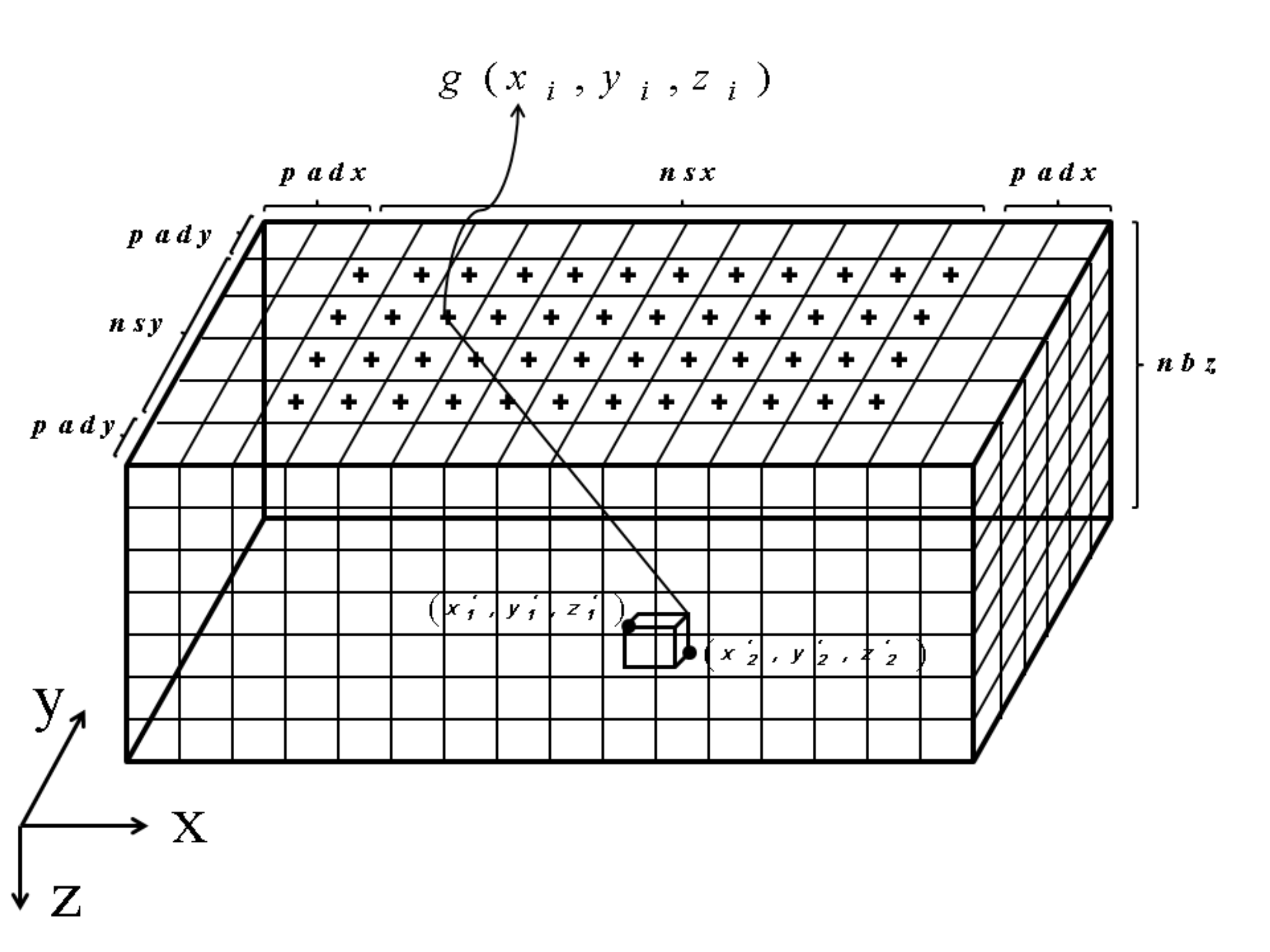}
\caption {\label{fig1} Discretization of the subsurface by rectangular prisms. $nsx$, and $nsy$ denote the number of gravity stations in the $x$ and $y$ directions, while  $nbz$ is the number of blocks in the (depth) $z$ direction. $padx$ and $pady$ denote the numbers of cells which may added around the gravity data grid in $x$ and $y$ directions, respectively.}
\end{figure}

Rectangular grid cells are commonly used for 3-D modelling of gravity sources. The subsurface under the survey area is divided into prisms of known sizes and positions. The unknown density contrasts within each prism define  the parameters to be estimated. Fig.~\ref{fig1} illustrates the  discretization of the subsurface by rectangular prisms. Gravity stations are located at the centers of the upper faces of  the prisms in the top layer. The cells are of equal size in each dimension,   $\Delta x =\Delta y=\Delta z$ where  $\Delta \cdot $ is the  distance  between gravity stations.  Extra cells may be added around the gravity data grid to reduce possible distortions in the reconstruction along the boundary \cite{Bou:2001}.

The vertical component of the gravitational attraction $g_i$ of a prism at point
$(x_i,y_i,z_i)$ is given by, \cite{Bou:2001}
\begin{eqnarray}\label{prism}
\frac{g_i}{\rho_j}&=&-\Gamma\sum_{p=1}^2\sum_{l=1}^2\sum_{s=1}^2\mu_{pls}\left[a_p\ln
\left(b_l+r_{pls}\right)+b_l\ln\left(a_p+r_{pls}\right)-c_s\arctan\left(\frac{a_pb_l}{c_sr_{pls}}\right)\right],
\end{eqnarray}
with
\begin{eqnarray}\mu_{pls}&=&\left(-1\right)^p\left(-1\right)^l\left(-1\right)^s \quad r_{pls}=\left(a_p^2+b_l^2+c_s^2\right)^\frac{1}{2}  \quad \mathrm{and}\\ 
a_p&=&x_i-x'_p,\quad  b_l=y_i-y'_l,\quad c_s=z_i-z'_s,\quad p,l,s=1,2.\end{eqnarray}
The coordinates of the eight corners for prism $j$ are denoted by  $(x'_p, y'_l, z'_s)$.
In \eqref{prism}  $\Gamma$ is the universal gravitational constant, $\rho_j$   is the density of the $j$th prism and $r_{pls}$ is the distance between one corner of the prism and the observation point.                                                                                             
The term on the right-hand side of  \eqref{prism}, which quantifies the contribution to the $i$th datum of  unit density in the $j$th cell, is denoted by the kernel  weight $G_{ij}$, and  is valid only at station $i$  for  cell $j$. The total response  for  station $i$  is obtained   by summing over all cells giving 
\begin{eqnarray}
g_i &=& \sum_{j=1}^{n} G_{ij} \rho_j, \quad i=1,\dots, m, 
\end{eqnarray}
leading to the linear equation
\begin{eqnarray}\label{matrix1}
\bfd = G \bfm, \quad m \ll n
\end{eqnarray}
Here we  use the standard notation that vector $\bfd \in\mathcal{R}^{m}$ is the set of measurements given by the $g_i$, and  $\bfm \in\mathcal{R}^{n}$ is the vector of unknown model parameters. 

Practical geophysical data are always contaminated by noise. Suppose that $\bfe \in \mathcal{R}^m$ represents the error in the measurements, assumed to be Gaussian and uncorrelated, then \eqref{matrix1} is replaced by 
\begin{eqnarray}\label{matrix2}
\bfdo = G \bfm +\bfe.
\end{eqnarray}
 The purpose of the gravity inverse problem is to find a geologically plausible density model $\bfm$ that reproduces $\bfd_{\mathrm{obs}}$ at the noise level.

\section{Focusing inversion methodology}\label{focusing}
An approximate solution for the ill-posed inverse problem  described by \eqref{matrix2} can be obtained by minimizing the  penalized least squares Tikhonov  functional defined by 
\begin{eqnarray}\label{tikhonov1}
\bfm(\alpha):=\argmin{\bfm}{\lbrace\|\Wd(G \bfm-\bfdo)\|_2^2 + \alpha^2 \|D(\bfm-\bfma)\|_2^2\rbrace}.
\end{eqnarray}
Here  $ \|\Wd(G\bfm-\bfdo)\|_2^2$ is the weighted data fidelity and $\|D(\bfm-\bfma)\|_2^2$ is the regularization term.  Data weighting matrix is given by  $\Wd=\mathrm{diag}(1/\eta_1,\dots,1/\eta_{m})$, where $\eta_{i}$ is the  standard deviation of the noise in the $i$th datum. $G\bfm$ is the vector of predicted data, $D$ is the regularization matrix and $\bfma$ is a given reference vector of a $\mathrm{priori}$ information for the model $\bfm$. In \eqref{tikhonov1} $\alpha$ is a regularization parameter which trades-off between the data fidelity and regularization terms. Introducing   $\tilde{G}:=\Wd G$  and  $\bfdto:=\Wd \bfdo$ in order to whiten the noise in the measurements $\bfdo$,  and shifting by the prior information through $ \bfmyalpha=\bfm-\bfma$, we find instead
\begin{eqnarray}\label{tikhonov2}
\bfmalpha:=\argmin{\bfy}\lbrace\|\tilde{G} \bfmyalpha-\tilde{\bfr}\|_2^2 + \alpha^2 \|D\bfmyalpha\|_2^2\rbrace, \quad \tilde{\bfr}=(\bfdto-\tilde{G}\bfma).
\end{eqnarray}
Under the assumption that the null spaces of $\tilde{G}$ and $D$ do not intersect,   $\bfm(\alpha)$  is
explicitly dependent on $\alpha$ and is defined in terms of the regularized inverse  $\tilde{G}(\alpha)$, 
\begin{eqnarray}
\bfmalpha  &=& (\tilde{G}^T\tilde{G}+\alpha^2D^TD)^{-1}\tilde{G}^T \tilde{\bfr}  = \tilde{G}(\alpha) \tilde{\bfr}, \quad \tilde{G}(\alpha):=(\tilde{G}^T\tilde{G}+\alpha^2D^TD)^{-1}\tilde{G}^T \label{reginv} \\
\bfm(\alpha)&=&\bfma+\bfmalpha = \bfma+ \tilde{G}(\alpha) \tilde{\bfr}. \label{tik2soln}
\end{eqnarray}
It is well-known that when the matrix $D$ is invertible the standard form transformation, \cite{Hansen}, yields the alternative but equivalent formulation 
\begin{eqnarray}
(\tilde G^T\tilde G + \alpha^2 D^TD) = D^T( (D^T)^{-1}\tilde G^T \tilde G D^{-1}+ \alpha^2 I_n)D.
\end{eqnarray}
The system describing the fidelity is replaced by the   right preconditioned matrix $\tilde{\tilde G}:=\tilde G D^{-1}$, giving the  regularized inverse $\tilde{\tilde G}(\alpha) := (\tilde{\tilde G}^T\tilde{\tilde G} + \alpha^2 I_n)^{-1} \tilde{\tilde G}^T $,  for which $ \bfz(\alpha)= D\bfy(\alpha)$ is defined by 
\begin{eqnarray}
\bfz(\alpha):=\argmin{\bfz}\lbrace\|\tilde{\tilde{G}} \bfz-\tilde{\bfr}\|_2^2 + \alpha^2 \|\bfz\|_2^2\rbrace.  \label{tikhonov3} \end{eqnarray}
Thus
\begin{eqnarray}
\bfm(\alpha)=\bfma + D^{-1}\bfz(\alpha).\label{tik3soln} 
\end{eqnarray}
Although analytically equivalent, numerical techniques to find \eqref{tik2soln} and \eqref{tik3soln} differ, for example using for \eqref{tik2soln} the generalized singular value decomposition, e.g. \cite{PaigeSau1}, for the matrix pair $[\tilde G, D]$, but the SVD of the $\tilde{\tilde G}$ for \eqref{tik3soln}, e.g.\cite{GoLo:96}. The solutions depend on the stability of these underlying decompositions, as well as the feasibility of calculating $D^{-1}$. 

Practically, the gravity inversion problem solves \eqref{tikhonov1} with an iteratively-defined operator, $D^{(k)}\in\mathcal{R}^{n \times n}$  given  by the product $D^{(k)}=W^{(k)}_{\mathrm{e}}W_{\mathrm{depth}}W_{\mathrm{hard}}$. While the depth weighting matrix \cite{Li:96}, $W_{\mathrm{depth}}=\mathrm{diag}(1/(z_{j})^\beta)$, and the hard constraint matrix, $W_{\mathrm{hard}}$ are independent of the iteration index,  the MS stabilizer matrix \cite{PoZh:99}, depends on the iteration. Specifically,  $W^{(k)}_{\mathrm{e}} =\mathrm{diag}\left((\bfm^{(k)}-\bfm^{(k-1)})^2+\epsilon^2 \right)^{-1/2}$,  $k>0$, with $W^{(0)}=I$ and $\bfm^{(0)}=\bfma$, see \cite{vatan:2014}.  The parameter $\epsilon>0$ is a focusing parameter which provides stability as $\bfm^{(k)}\rightarrow\bfm^{(k-1)}$ and parameter $\beta$  determines the weight on the cell $j$ with mean depth $z_{j}$. The hard constraint matrix $W_{\mathrm{hard}}$ is  initialized as the identity matrix, with    $(W_{\mathrm{hard}})_{jj} =H$ , where $H$ is a large number which then forces  $(\bfma)_j=\rho_j$  for those $j$ where  geological and geophysical information are able to provide the value of the density of cell $j$.   In order to recover  a feasible image of the subsurface lower and upper density bounds   $[\rho_{\mathrm{min}}, \rho_{\mathrm{max}}]$  are imposed. During the inversion process if a given density value falls outside the bounds, the value at that  cell  is projected back to the nearest  constraint value. Furthermore,  the algorithm terminates when the solution either reaches the noise level, i.e.  $\chi^2_{\text{Computed}}:=\| (\bfdo)_{i}-(\bfdp)_{i} / \eta_{i} \|^2_2 \leq m+\sqrt {2 m}$, or a maximum number of iterations is reached. 

The iterative formulation of \eqref{tikhonov3}, given $\{\alpha^{(k)}, k>0 \}$,  is now clear. We set regularizer $D^{(k)}=D(\bfm^{(k)},\bfm^{(k-1)})$ and $\tilde{\bfr}^{(k)} = \bfdo-\tilde{G}\bfm^{(k)}$ for $k>1$, initialized with $\tilde{\bfr}^{(0)}=\bfdo-\tilde{G}\bfma$ and $D^{(0)}=W_{\mathrm{depth}}$, yielding the regularization parameter dependent updates  
\begin{eqnarray}\label{normalk}
\bfz(\alpha^{(k+1)}) &=&( \tilde{\tilde G}^T\tilde{\tilde G}+(\alpha^{(k)})^2I_n)^{-1}\tilde{\tilde G} \tilde{\bfr}^{(k)}, \\
\bfm^{(k+1)}&=&\bfm^{(k)} + (D^{(k+1)})^{-1}\bfz(\alpha^{(k+1)}).\label{tikk} 
\end{eqnarray}
Using the  SVD for the matrix $\tilde{\tilde G}$, see \ref{appSVD}, \eqref{normalk} can be written as
\begin{eqnarray}\label{svdsolution}
\bfz(\alpha ^{(k+1)}) &=& \sum_{i=1}^{m} \frac{\sigma^2_i}{\sigma^2_i+(\alpha^{(k)})^2} \frac{\bfu^T_{i}\tilde{\bfr}^{(k)}}{\sigma_i} \bfv_i
\end{eqnarray}
This  formulation  \eqref{svdsolution} demonstrates that we may efficiently accomplish the solver through use of the SVD in place of the GSVD. 

Still,  the algorithm suggested by \eqref{normalk}-\eqref{tikk} also requires estimation of the parameter $\alpha^{(k)}$ which further complicates the solution process. First, an approach for determining or describing an optimal $\alpha$ must be adopted and rationalized. Second, regardless of the criterion that is chosen for finding $\alpha$, the implementation  requires calculating $\bfm(\alpha)$ for multiple choices of $\alpha$. It is therefore crucial to have   an effective criterion for defining an optimal $\alpha$ at each step.

\section{Regularization parameter estimation}\label{parameter estimation}
Effective and efficient regularization parameter estimation for Tikhonov regularization is well-described in the literature e.g. \cite{Hansen,Vogel:2002}. In the context of the gravity inversion problem the regularization parameter  $\alpha$ is required 
at each  iteration $k$, and thus the problem of finding the \textit{ optimal} $\alpha:=\alphao$ efficiently is even more crucial. One approach that has been previously adopted in the literature is an iterated Tikhonov procedure in which $\alpha^{(k)}$ is chosen to converge geometrically, e.g. $\alpha^{(k)}=\alpha^{(1)}q^{(k)}$ for a decreasing geometric sequence $q^{(k)}$, e.g. $q^{(k)}=2^{-k}$, \cite{Tikhonov,Zhdanov}, hence eliminating the need to estimate the parameter for other than the first step.  Our results will show that this would not be useful here.  Assuming then that $\alpha$ is updated each step,  the most often used method  for potential field data inversion is the MDP. Yet it is well-known that the MDP always leads to an over estimation of the regularization parameter,  e.g. \cite{kilmer:2001}, and hence an over smoothing of the solution. Further,  the LC and GCV are techniques which extend easily for underdetermined systems, without any additional analysis, and were therefore considered in \cite{vatan:2014}. On the other hand, the UPRE and $\chi^2$ techniques were developed for the solution of underdetermined problems, extending prior results for consistent or overdetermined systems, and carefully validated for their use in 2D focusing inversion \cite{VRA2014}. These results  indicate a preference for the UPRE and $\chi^2$ techniques. Thus here we focus on the comparison of the established MDP with the UPRE and $\chi^2$ techniques for  3D  potential field data inversion. Because the UPRE and $\chi^2$ techniques are less well-known for this problem domain, we briefly describe the rationale for the UPRE and $\chi^2$ techniques, but leave the presentation of the formulae to \ref{regparam} and point to \cite{VRA2014} for the derivations.  We note that as with the MDP, it is assumed that an estimate of the noise level in the data is provided.

\subsection{Unbiased predictive risk estimator}\label{UPRE}
Noting that the optimal $\alphao$ should minimize the  error between the Tikhonov regularized solution $ \bfz(\alpha) $ and the exact solution $\bfze$,  the purpose is to develop a method for effectively estimating this optimal $\alpha$ without knowledge of $\bfze$ through use of the  measurable residual and the statistical estimator of the mean squared norm of the error,   \cite{Vogel:2002}. Specifically, with $H(\alpha)=\tilde{\tilde G}\tilde{\tilde G}(\alpha)$, the 
predictive error $\bfp(\bfz(\alpha))$  given by 
\begin{eqnarray}\label{prederror}
\bfp(\bfz(\alpha)) &:= \tilde{\tilde G}\bfz(\alpha)-\bfrte=\tilde{\tilde G}\tilde{\tilde G}(\alpha)\tilde{\bfr}- \bfrte=
 (H(\alpha) -I_m)\bfrte +H(\alpha) \tilde{\bfe}, 
\end{eqnarray}
is not available, 
 but the  residual 
\begin{eqnarray}\label{measurable residual}
R(\bfz(\alpha)):= \tilde{\tilde{G} }\bfz(\alpha) - \tilde{\bfr} = (H(\alpha)-I_m)\tilde{\bfr} = (H(\alpha)-I_m)(\bfrte+\tilde{\bfe}),
\end{eqnarray}
is measurable. Thus an 
estimate of the 
 mean squared norm 
\begin{eqnarray}\label{predrisk}
\frac {1}{m}\| \bfp(\bfz(\alpha))  \|_2^2 &=& \frac {1}{m}\| (H(\alpha) -I_m)\bfrte +H(\alpha) \tilde{\bfe}  \|_2^2, 
\end{eqnarray}
is obtained via the mean squared norm for $R(\bfz(\alpha))$ and some algebra that employs the Trace Lemma \cite{Vogel:2002}. Then, the  optimal regularization parameter is selected such that
\begin{eqnarray}\label{opt1}
\alphao=\argmin{\alpha}\lbrace \frac {1}{m}\| \bfp(\bfz(\alpha))  \|_2^2  \rbrace= \argmin{\alpha}\lbrace U(\alpha)  \rbrace,
\end{eqnarray}
where 
\begin{eqnarray}\label{upre}
U(\alpha)=\|\tilde{\tilde{G}} \bfz(\alpha) - \tilde{\bfr}\|_2^2 +2\mathrm{trace}(H(\alpha))-m.
\end{eqnarray}
is the functional to be minimized for the UPRE technique to find $\alphao$. This functional can be evaluated in terms of the SVD, as indicated in \eqref{upresvd}.

\subsection{$ \chi^2 $ principle}\label{chi2}
The $\chi^2$ principle  is a generalization of the MDP. Whereas the MDP is obtained under the assumption that $\alphao$ should yield a fidelity term that follows a $\chi^2$ distribution with $m-n$ degrees of freedom, for overdetermined systems, the $\chi^2$ principle for regularization parameter estimation considers the entire Tikhonov functional. For weighting of the data fidelity by a known Guassian noise distribution on the measured data and, when the stabilizing term is considered to be weighted by unknown inverse covariance information on the model parameters, the minimum of the Tikhonov functional becomes a random variable that follows a $\chi^2$-distribution with $m$ degrees of freedom, \cite{mere:09,VRA2014}, a result that holds also for underdetermined systems, which is not the case for the MDP. Specifically for the MDP one seeks in general 
\begin{eqnarray}\label{residual}
\|\tilde{\tilde{G}}\bfz(\alpha)-\tilde{\bfr}\|_2^2= m-n, \quad m\ge n, 
\end{eqnarray}
which is then usually replaced by an estimate based on the variance when $m<n$, see e.g.  \cite{Far:2004},  while for the $\chi^2$ principle we seek 
\begin{eqnarray}\label{tikhonov4}
P(\bfm(\alpha))=\|\Wd( G\bfm-\bfdo)\|_2^2 + \alpha^2\| D(\bfm-\bfma)   \|^2 = m,
\end{eqnarray}
which is under the assumption that $\alpha^2I$ effectively whitens the noise in the estimate for $\bfm$ around the mean $\bfma$. These yield the formulae
\eqref{MDPSVD}  and  \eqref{chi2GSVD}  for the MDP and $\chi^2$ principle, respectively, when used with the SVD.

\section{Synthetic examples}\label{synthetic}

\subsection{Synthetic example: Dike}\label{dike}
\begin{figure}
\subfigure{\label{2a}\includegraphics[width=.49\textwidth]{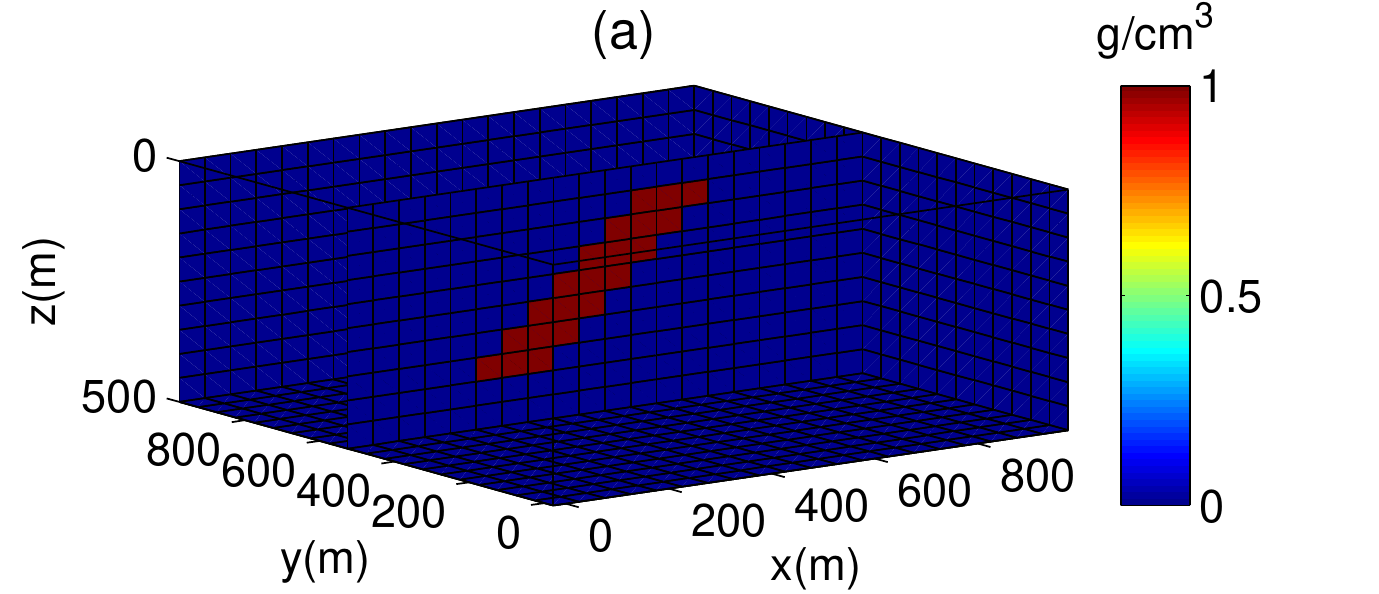}}
\subfigure{\label{2b}\includegraphics[width=.49\textwidth]{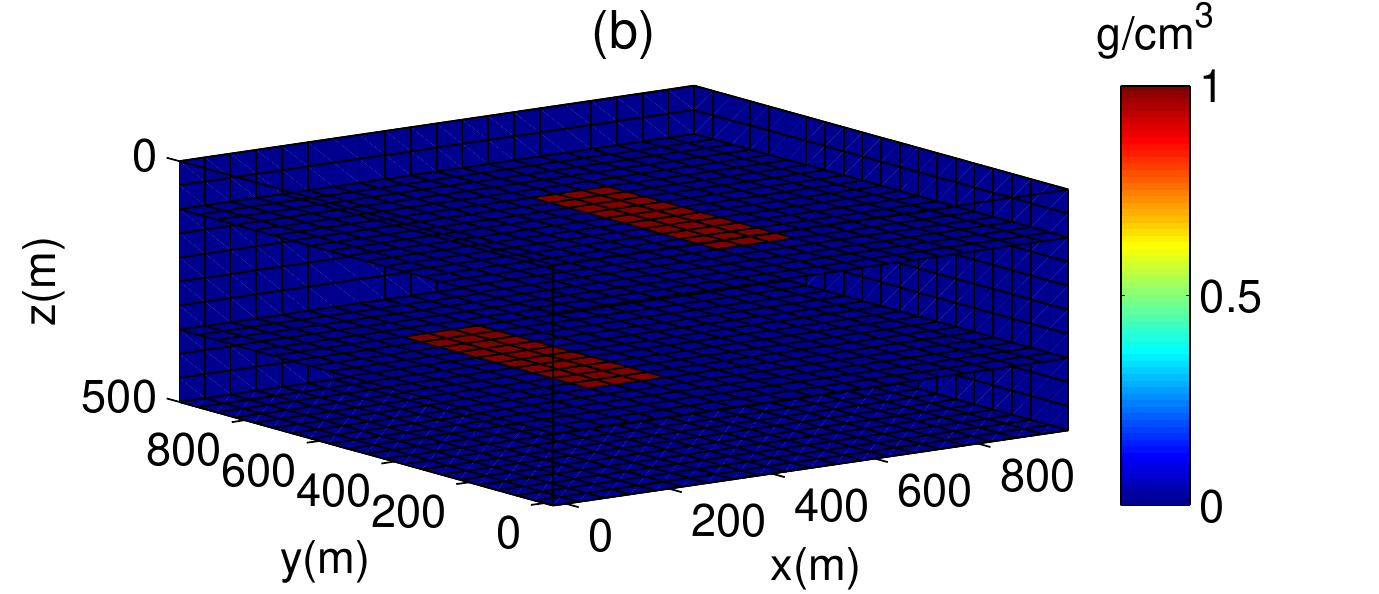}}
\caption {\label{fig2} Model of a dipping dike on  an homogeneous background. Fig.~\ref{2a}: cross-section at $y=525$~m;  Fig.~\ref{2b}: plane-sections at $z=100$~m and $z=350$~m. The density contrast of the dike is $1$~g$/$cm$^3$.}
\end{figure}

\begin{figure}
\subfigure{\label{3a}\includegraphics[width=.4\textwidth]{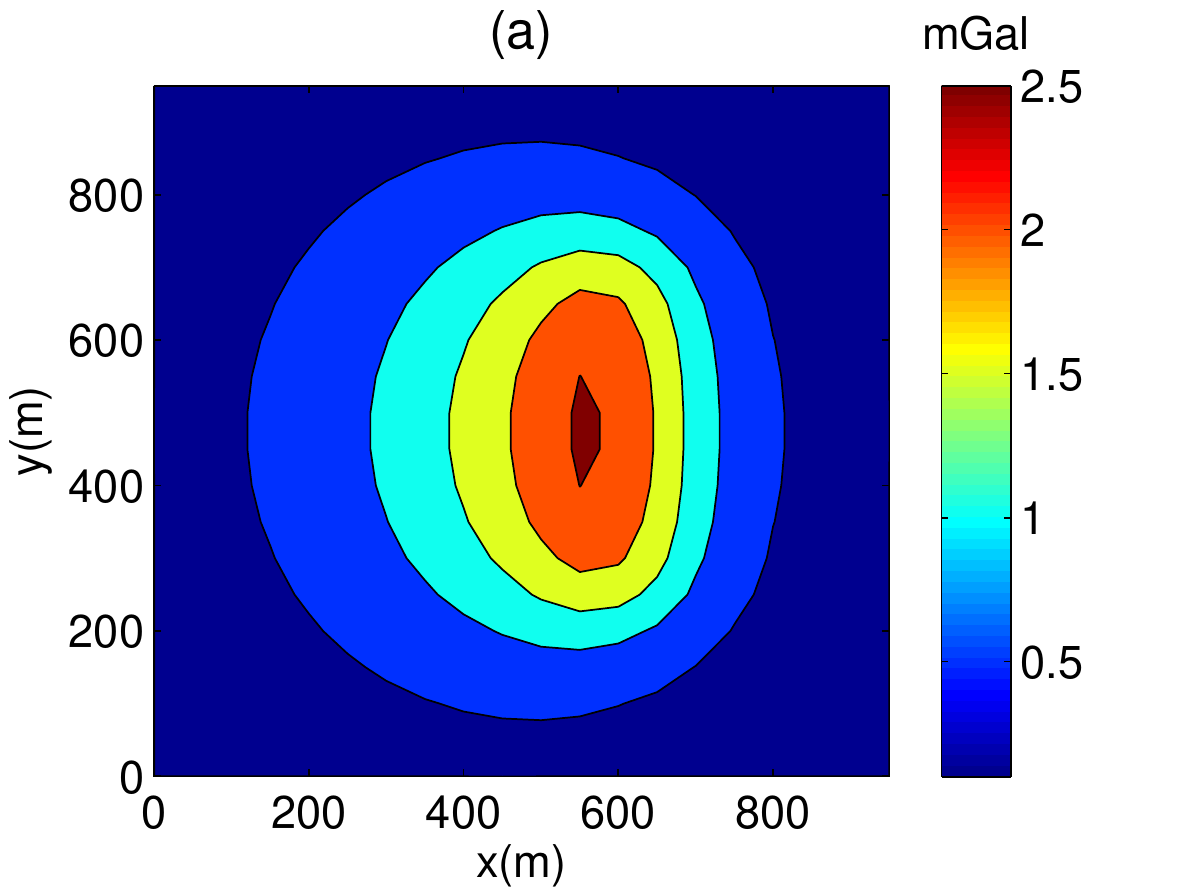}}
\subfigure{\label{3b}\includegraphics[width=.4\textwidth]{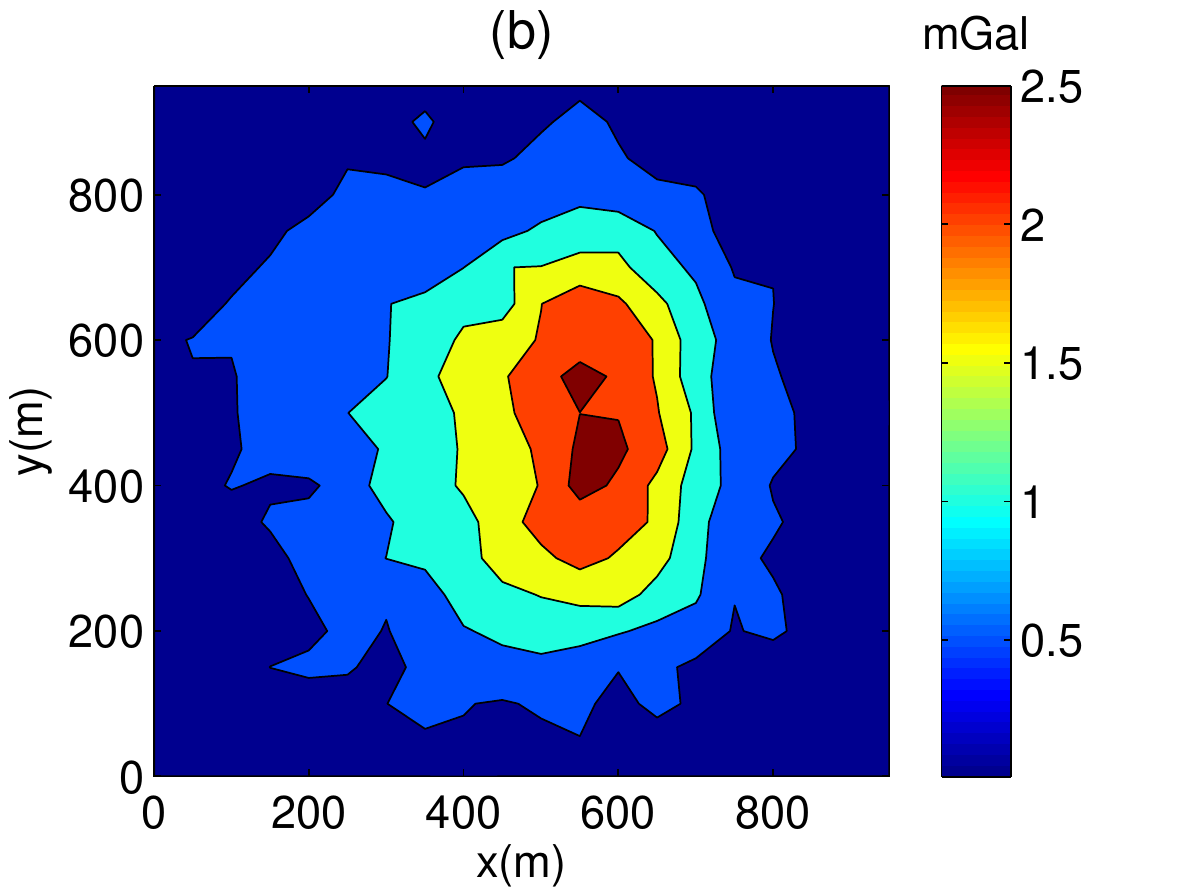}}
\caption {\label{fig3} Anomaly due to the dike model shown in Fig.~\ref{fig2}. Fig.~\ref{3a}: noise free data; Fig.~\ref{3b}:  data with added noise for $(\eta_1, \eta_2)=(0.02, 0.005)$.}
\end{figure}

The first model which is used for testing the reliability of the introduced parameter-choice methods is the dipping dike. Figs~\ref{2a}-\ref{2b} show  the cross  and plane sections of this model. It has density  contrast $1$~g$/$cm$^3$ on  an homogeneous background. Simulation data, $\bfd$, are calculated over a $20$ by $20$ grid with $\Delta =50$~m  on the surface, Fig.~\ref{3a}. In generating noise-contaminated data we generate a random matrix $\Theta$ of size $m \times 10$ using the MATLAB function $\bf{randn}$. Then setting $\bfdo^c=\bfd+(\eta_1(\bfd)_i+\eta_2\| \bfd\|)\Theta^c$, $c=1:10,$ generates $10$ copies of the right-hand side vector. The inversion results are presented for  $3$ noise realizations, namely $(\eta_1,  \eta_2)=(0.01, 0.001)$; $(\eta_1, \eta_2)=(0.02, 0.005)$; and  $(\eta_1,\eta_2)=(0.03, 0.01)$.  Fig.~\ref{3b} shows an example of noise-contaminated data for one right-hand side,  here $c=4$, for the second noise realization.

For inversion the subsurface is divided into $ 20 \times 20 \times 10 = 4000$ cells  each with $\Delta= 50$~m. The iterations are initialized with  $\bfma=\bf{0}$ and $W_{\mathrm{e}}= W_{\mathrm{hard}}= I_n$. Realistic bounds on the density are imposed by choosing $\rho_{\mathrm{min}}= 0$~g$/$cm$^3$ and $ \rho_{\mathrm{max}}=1$~g$/$cm$^3$. For all inversions the coefficient $\beta$ in $W_{\mathrm{depth}}$ and the focusing parameter $\epsilon$ are fixed at $0.8$ and $0.02$, respectively. The algorithm terminates when  $\chi^2_{\mathrm{Computed}} \leq 429$ or a maximum number of iterations, $K$, is reached. Here $K = 100$. The inversion is performed for all noise realization choices given by the $(\eta_1, \eta_2)$ pairs, and all $10$ random copies of the noise simulation in each case. The  following average values are calculated for all $10$ simulations in each case:  (i) the average regularization parameter at the final value, $\alpha^{(K)}$,  (ii)  the average number of iterations $K$  required for convergence, and (iii) the average relative error of the  reconstructed model,$\| \bfme-\bfm^{(K)}\|_2 /\| \bfme\|_2 $. The results are presented in Tables~\ref{tab1} - \ref{tab3}, for parameter estimation using the $\chi^2$ principle, the UPRE method, and the MDP method, respectively. Frequently, in potential field data inversion,  the initial value of the regularization parameter is taken to be large \cite{Far:2004}, i.e. at the first step no parameter choice method is required. We consistently initialize $\alpha^{(1)}$ for all methods  using the already known  singular values of the matrix $\tilde{\tilde{G}}$.  Specifically we take $\alpha^{(1)}=(n/m)^{\gamma}(\mathrm{max}(\sigma_i)/\mathrm{mean}(\sigma_i))$. Our investigations   show that $\gamma$ can be chosen such that $0 \leq \gamma \le 2$.

\begin{table}
\caption{The inversion results obtained by inverting the data from the dike contaminated with the  first noise level, $(\eta_1,  \eta_2)=(0.01, 0.001)$, average(standard deviation) over $10$ runs.}\label{tab1}
\begin{tabular}{c  c  c  c  c }
\hline
Method&     $\alpha^{(1)},\gamma=1.5$& $\alpha^{(K)}$& Relative error& Number of iterations \\ \hline
$\chi^2$ principle&  4737& 287(4.3)& 0.7752(0.0048)& 80.8(6.6)\\
UPRE&  4737& 63(0.001)& 0.7699(0.0050)& 58.9(4.8)\\
MDP&  4737& 215(8.4)& 0.7731(0.0051)& 100\\ \hline
\end{tabular}
\end{table}

\begin{table}
\caption{The inversion results obtained by inverting the data from the dike contaminated with the  second noise level, $(\eta_1, \eta_2)=(0.02, 0.005)$, average(standard deviation)  over $10$ runs..}\label{tab2}
\begin{tabular}{c  c  c  c  c }
\hline
Method&     $\alpha^{(1)},\gamma=1.5$& $\alpha^{(K)}$& Relative error& Number of iterations \\ \hline
$\chi^2$ principle&  4847& 66(6.7)& 0.7672(0.0089)& 6.2(0.9)\\
UPRE&  4847& 17.6(1.0)& 0.7662(0.0086)& 6.6(0.7)\\
MDP&  4847& 47.1(2.9)& 0.7808(0.0107)& 12.7(2.6)\\ \hline
\end{tabular}
\end{table}

\begin{table}
\caption{The inversion results obtained by inverting the data from the dike contaminated with the  third noise level, $(\eta_1,  \eta_2)=(0.03, 0.01)$, average(standard deviation)  over $10$ runs.}\label{tab3}
\begin{tabular}{c  c  c  c  c }
\hline
Method&     $\alpha^{(1)},\gamma=1.5$& $\alpha^{(K)}$& Relative error& Number of iterations \\ \hline
$\chi^2$ principle&  4886& 40.8(5.5)& 0.7574(0.0132)& 3\\
UPRE&  4886& 15.8(6.8)& 0.7404(0.0149)& 3.1(0.31)\\
MDP&  4886& 36.6(12.2)& 0.7786(0.0133)& 3.1(0.31)\\ \hline
\end{tabular}
\end{table}

The results in Tables~\ref{tab1}-\ref{tab3} show that both the $\chi^2$ and MDP methods lead to an overestimate of  the regularization parameter as compared to that obtained with the UPRE. On the other hand, with respect to the relative error of the reconstructed model, both the $\chi^2$ and  UPRE methods lead to reduced error as compared to the MDP. Furthermore, they both require  fewer iterations as compared to the MDP and the cost per iteration for the $\chi^2$ method is   cheaper than that for the  UPRE, requiring  just an efficient  root-finding algorithm while the UPRE relies on an  estimate of  $U(\alpha)$ on a range of $\alpha$.

\begin{figure}
\subfigure{\label{4a}\includegraphics[width=.49\textwidth]{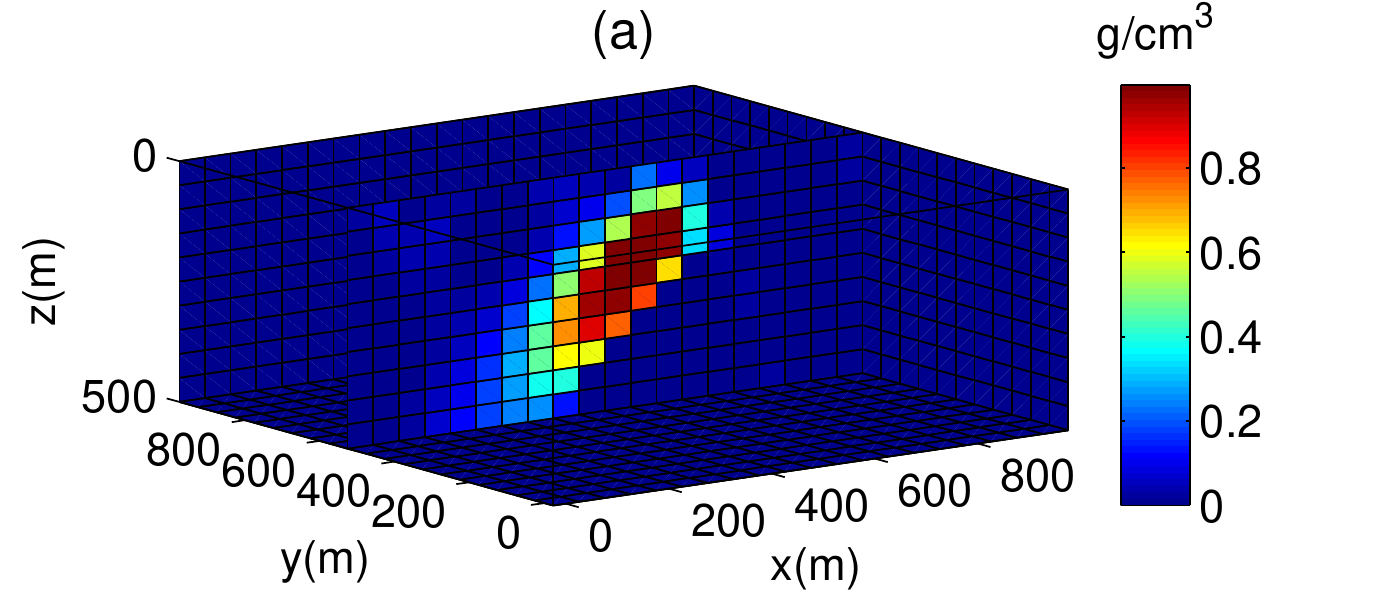}}
\subfigure{\label{4b}\includegraphics[width=.49\textwidth]{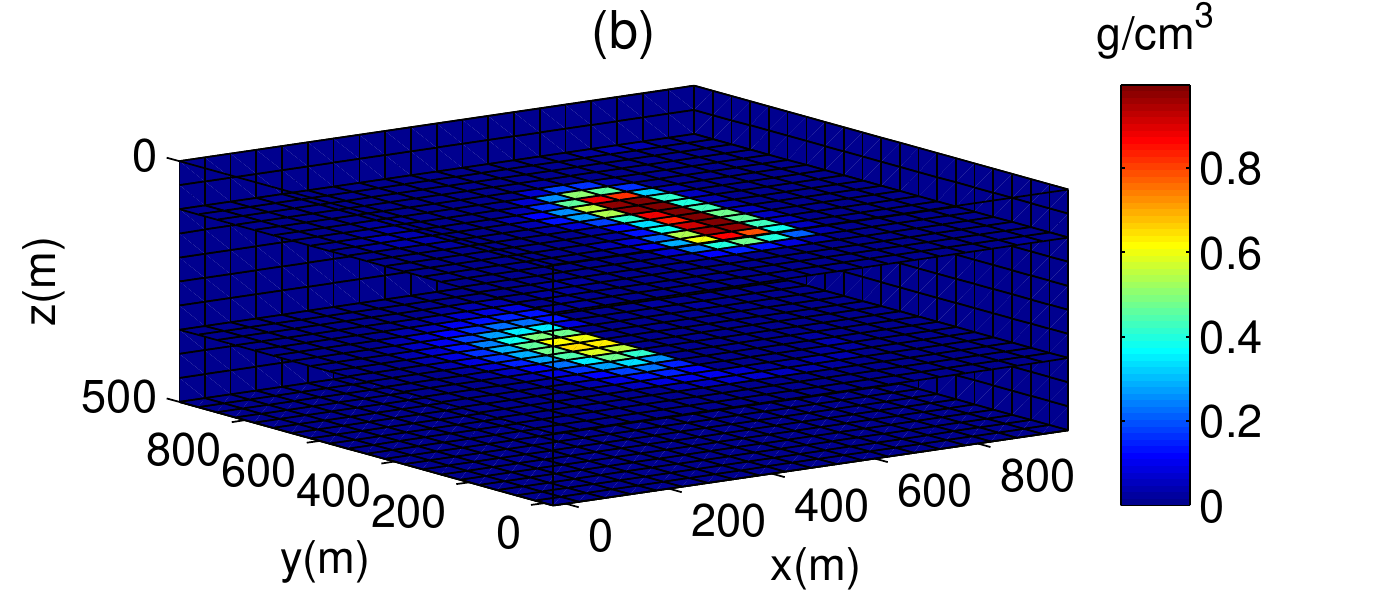}}
\subfigure{\label{5a}\includegraphics[width=.49\textwidth]{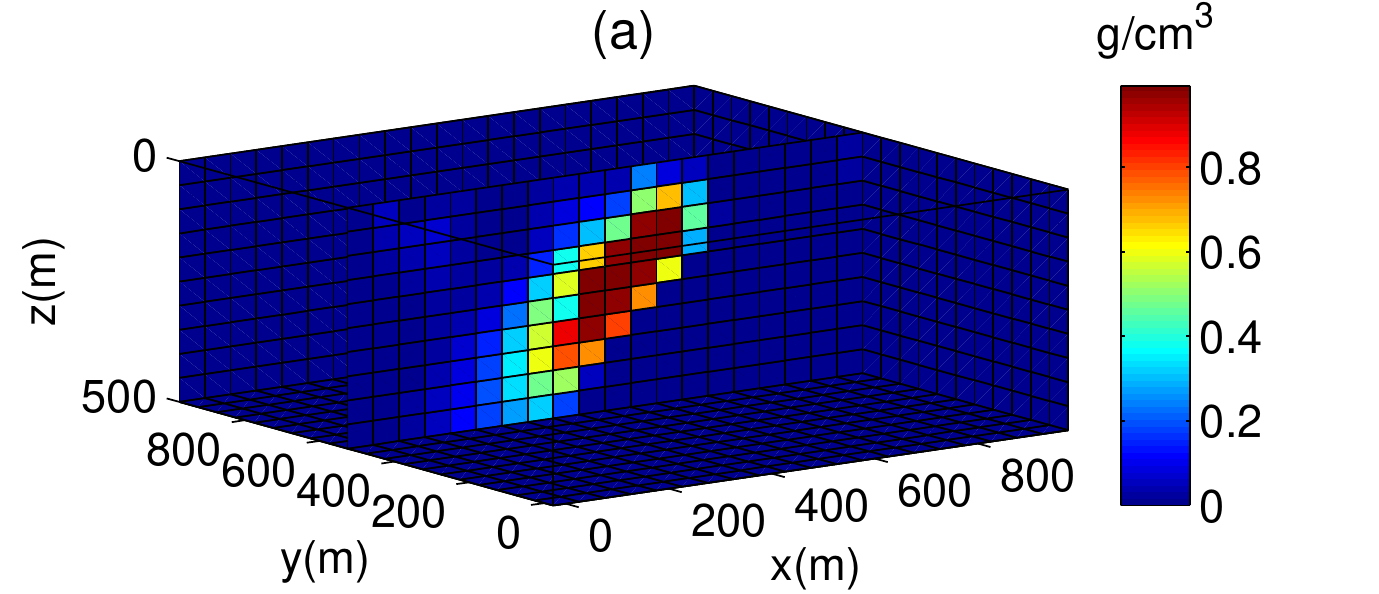}}
\subfigure{\label{5b}\includegraphics[width=.49\textwidth]{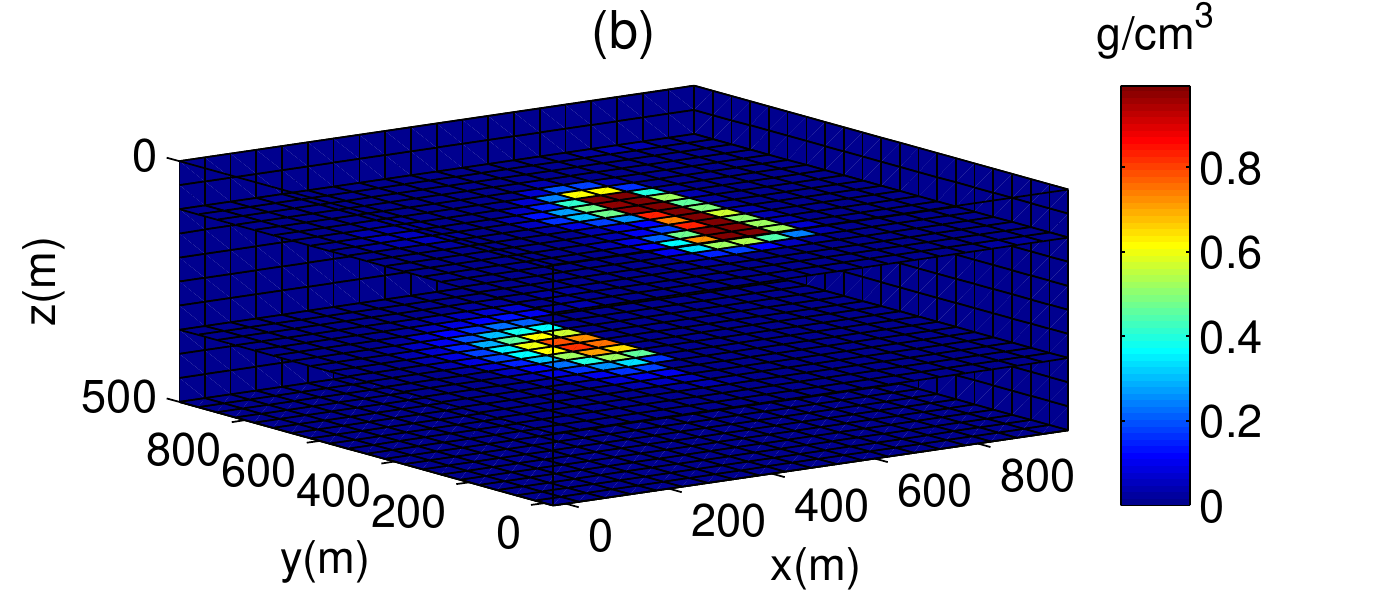}}
\subfigure{\label{6a}\includegraphics[width=.49\textwidth]{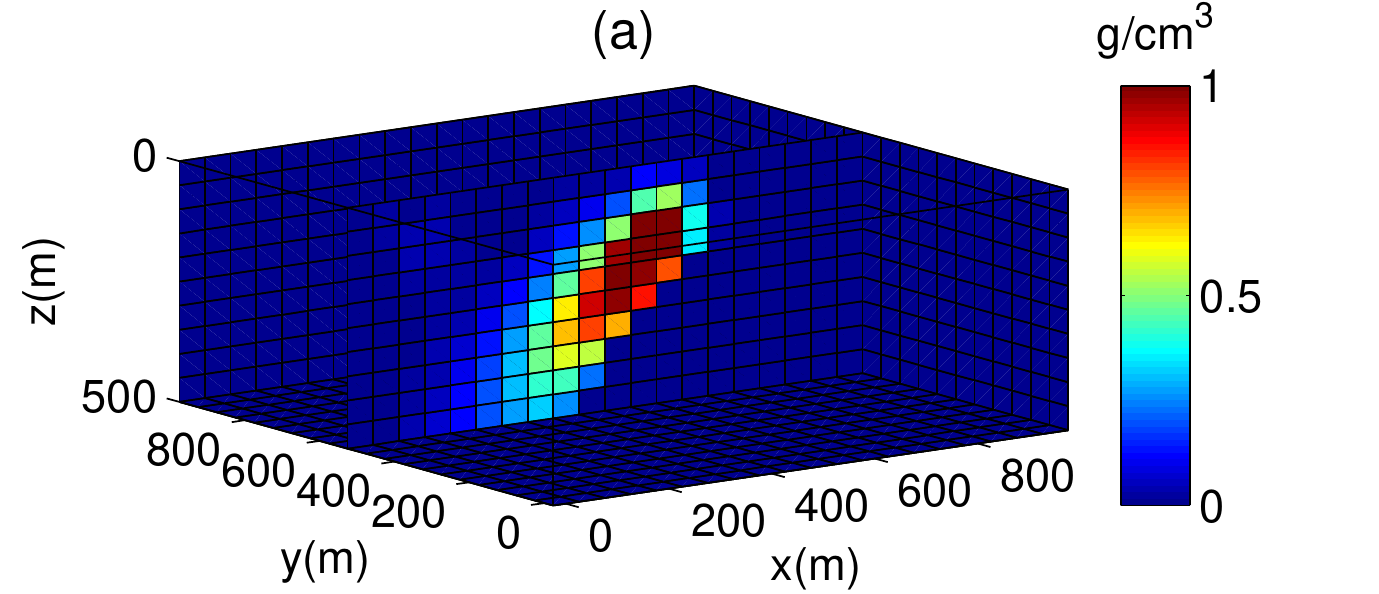}}
\subfigure{\label{6b}\includegraphics[width=.49\textwidth]{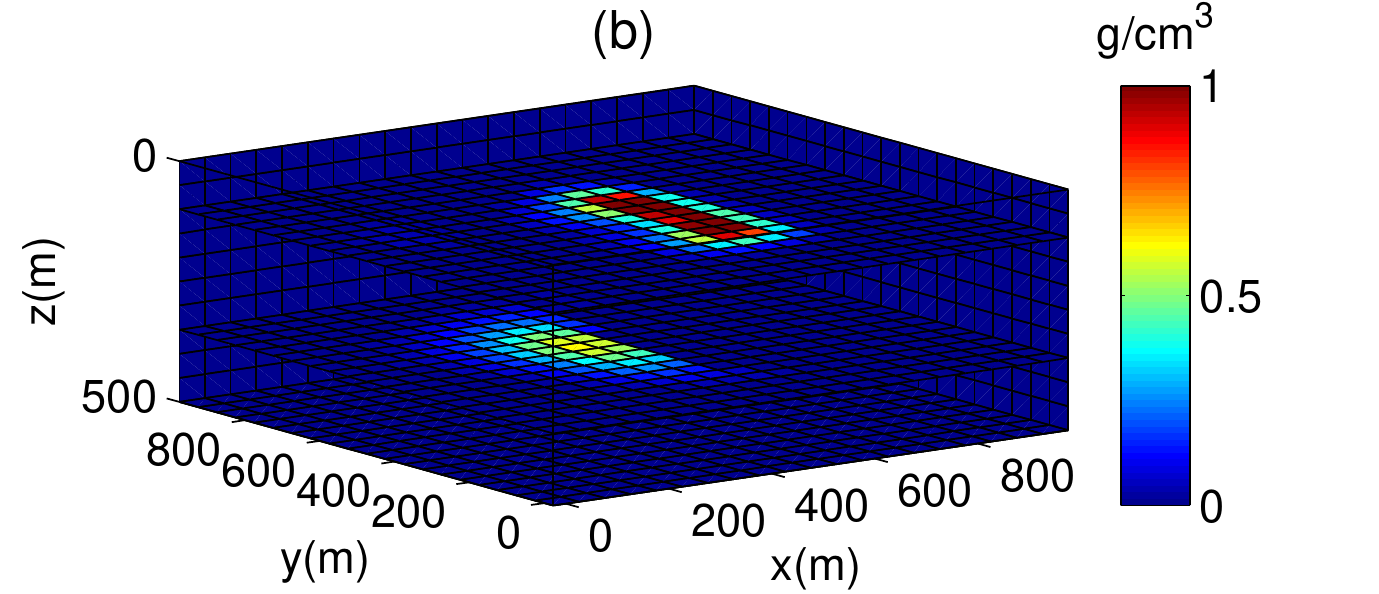}}
\caption {\label{fig4} The results obtained by inverting the data shown in Fig.~\ref{3b} using  the $\chi^2$ principle, the UPRE and the MDP  as the parameter-choice method, respectively. Figs~\ref{4a}, \ref{5a}, \ref{6a}: the cross-section at $y=525$~m in each case, respectively and in  Figs~\ref{4b}, \ref{5b}, \ref{6b}: the  plane-sections at $z=100$~m and $z=350$~m for the same cases.}
\end{figure}
\begin{figure}
\subfigure{\label{4c}\includegraphics[width=.35\textwidth]{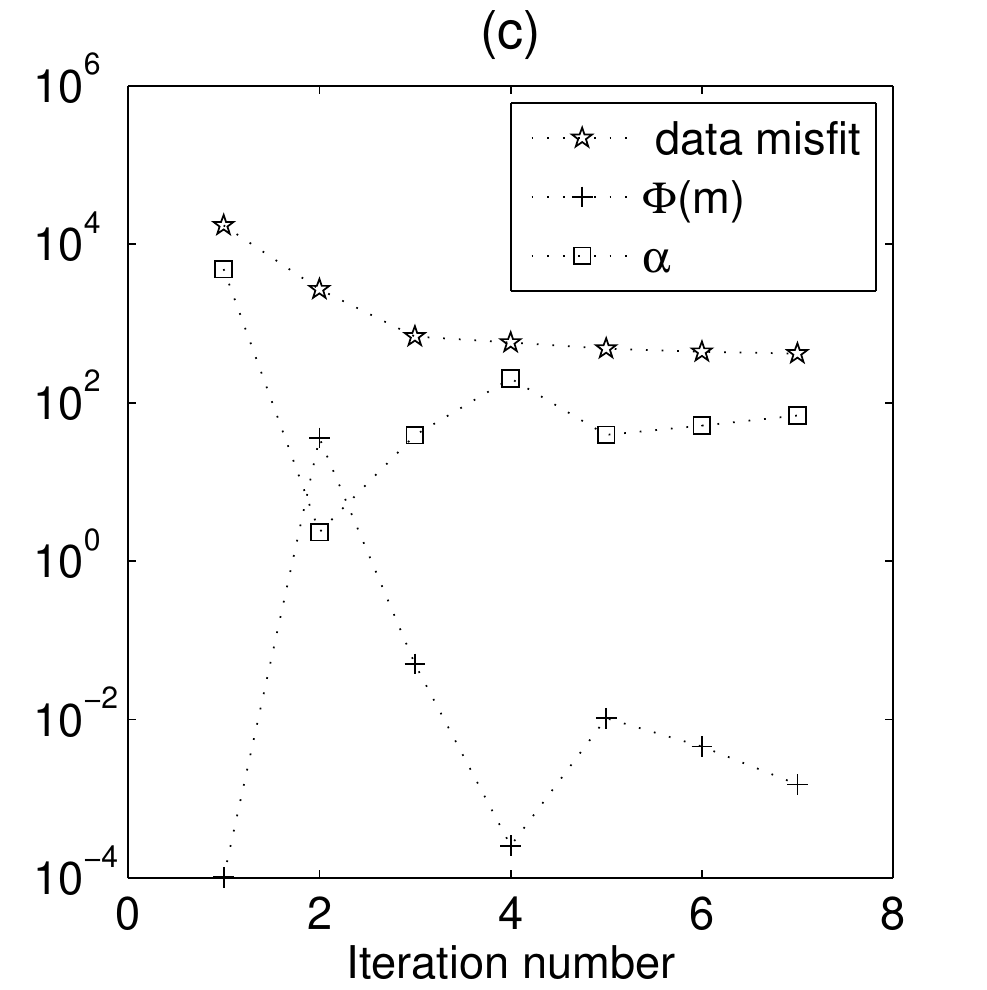}}
\subfigure{\label{4d}\includegraphics[width=.35\textwidth]{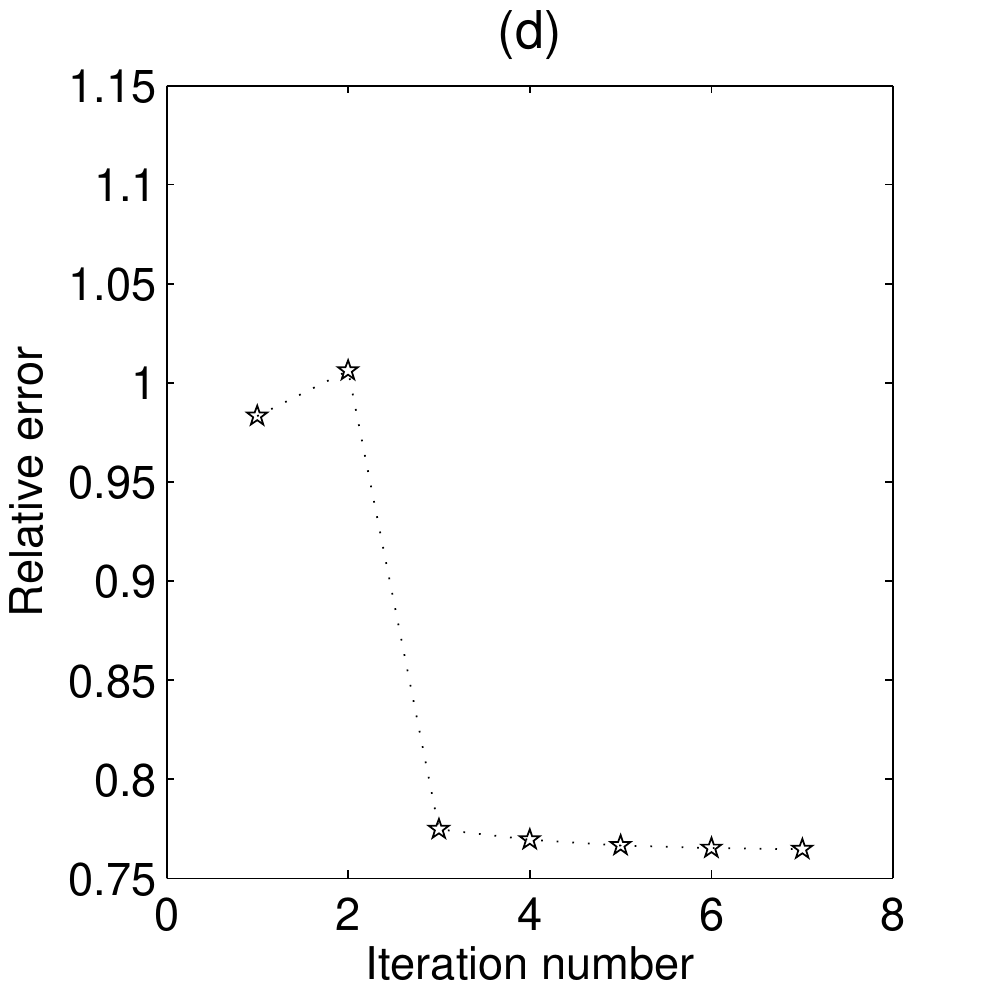}}\\
\subfigure{\label{5c}\includegraphics[width=.35\textwidth]{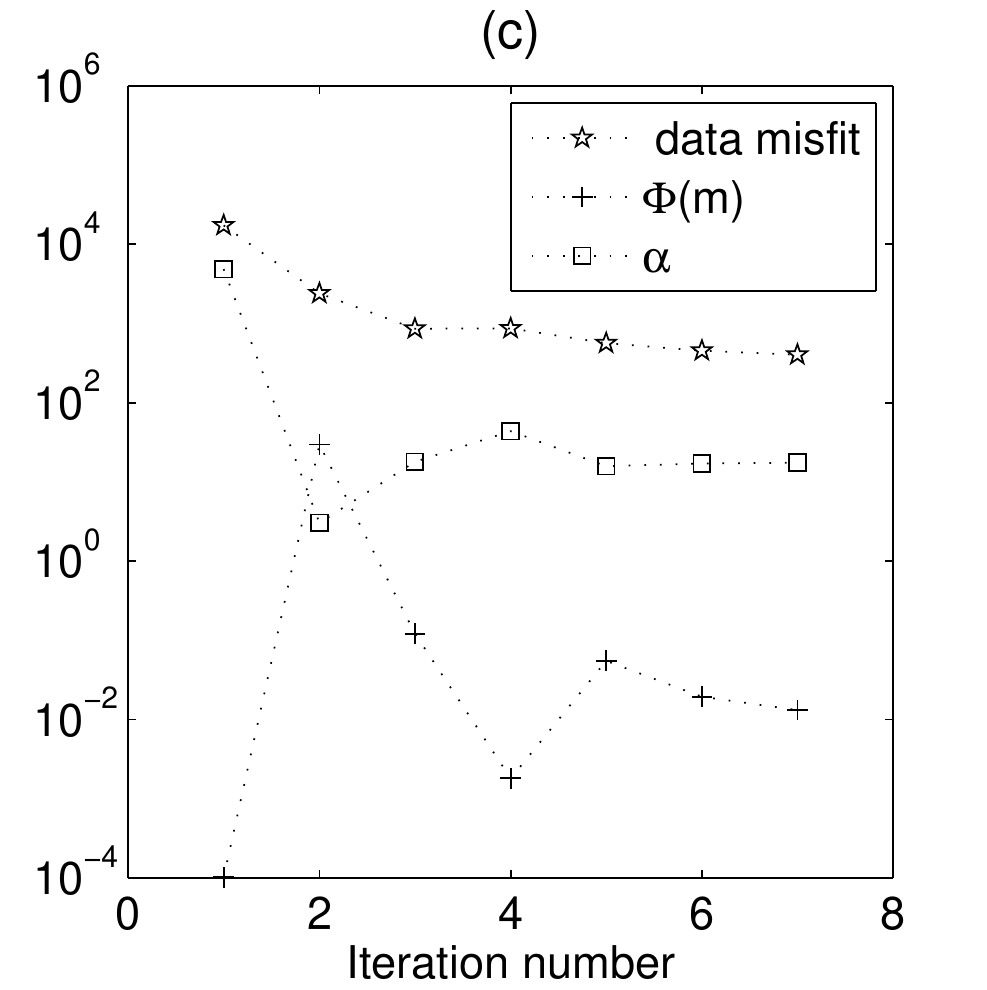}}
\subfigure{\label{5d}\includegraphics[width=.35\textwidth]{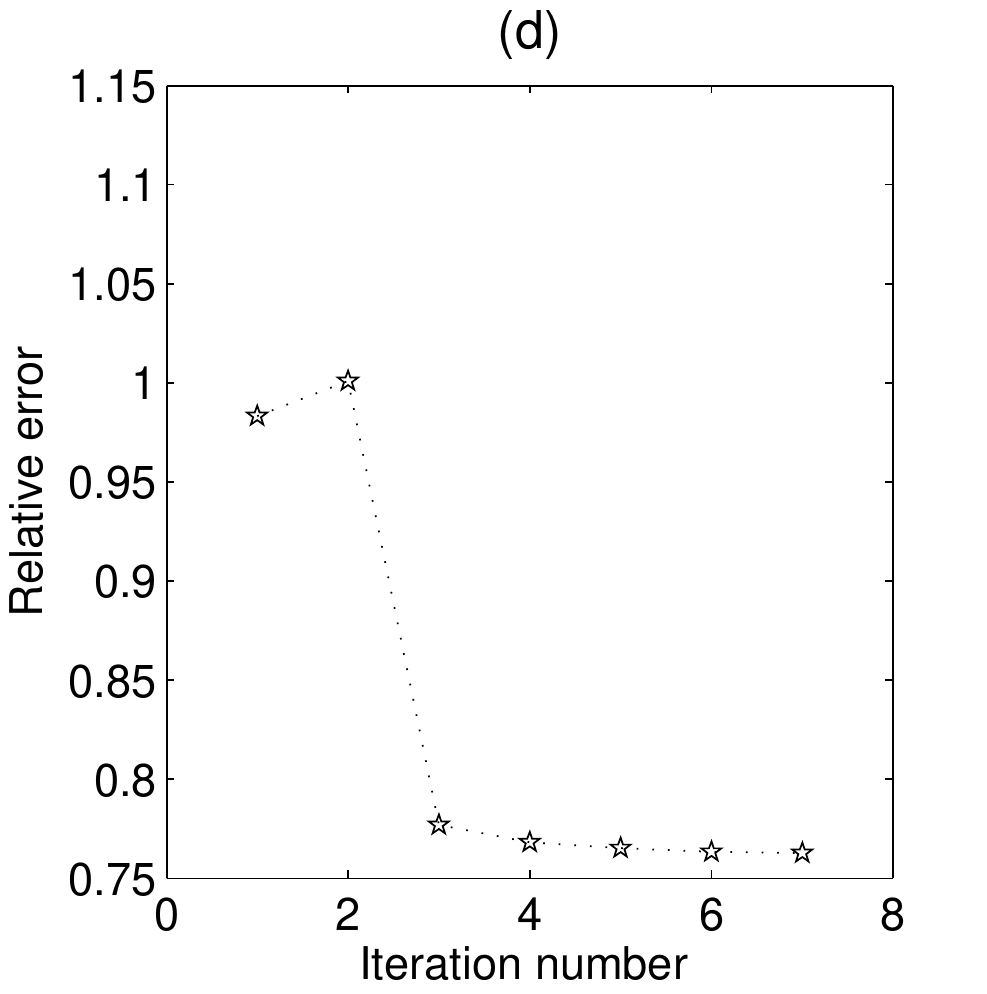}}\\
\subfigure{\label{6c}\includegraphics[width=.35\textwidth]{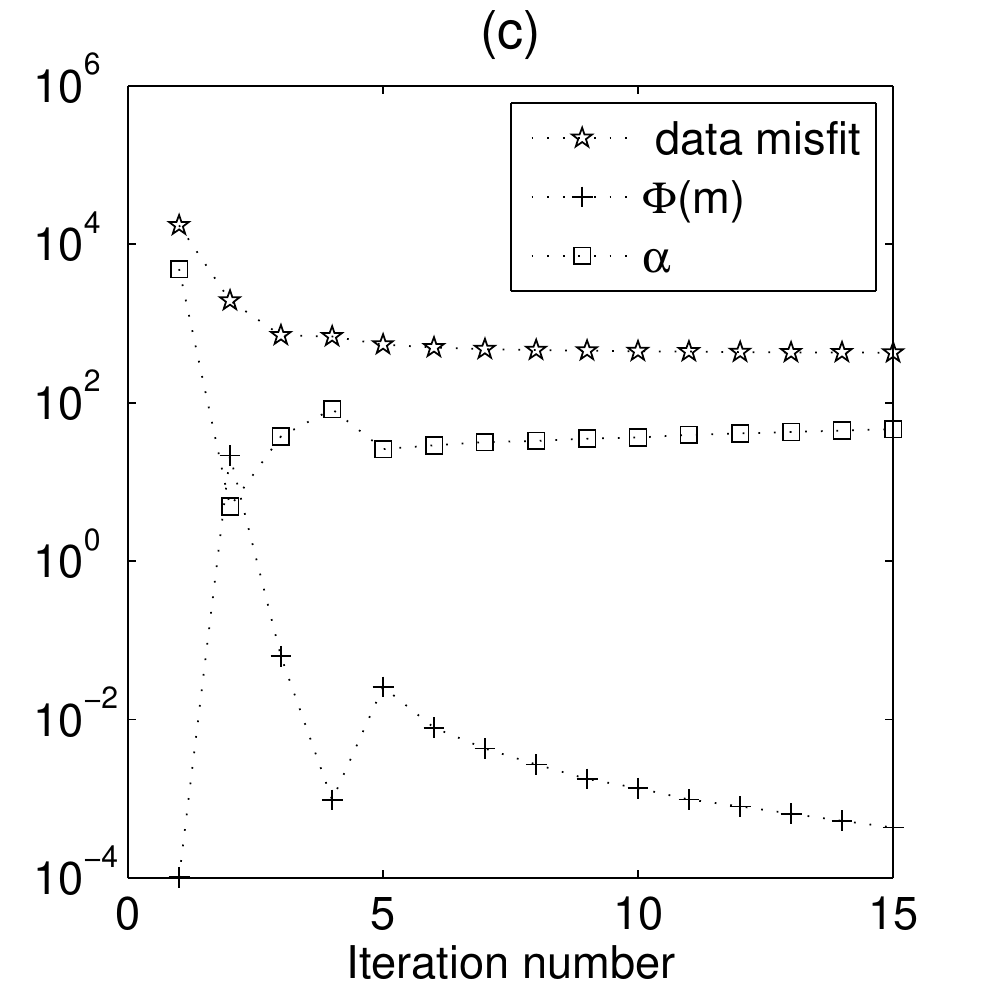}}
\subfigure{\label{6d}\includegraphics[width=.35\textwidth]{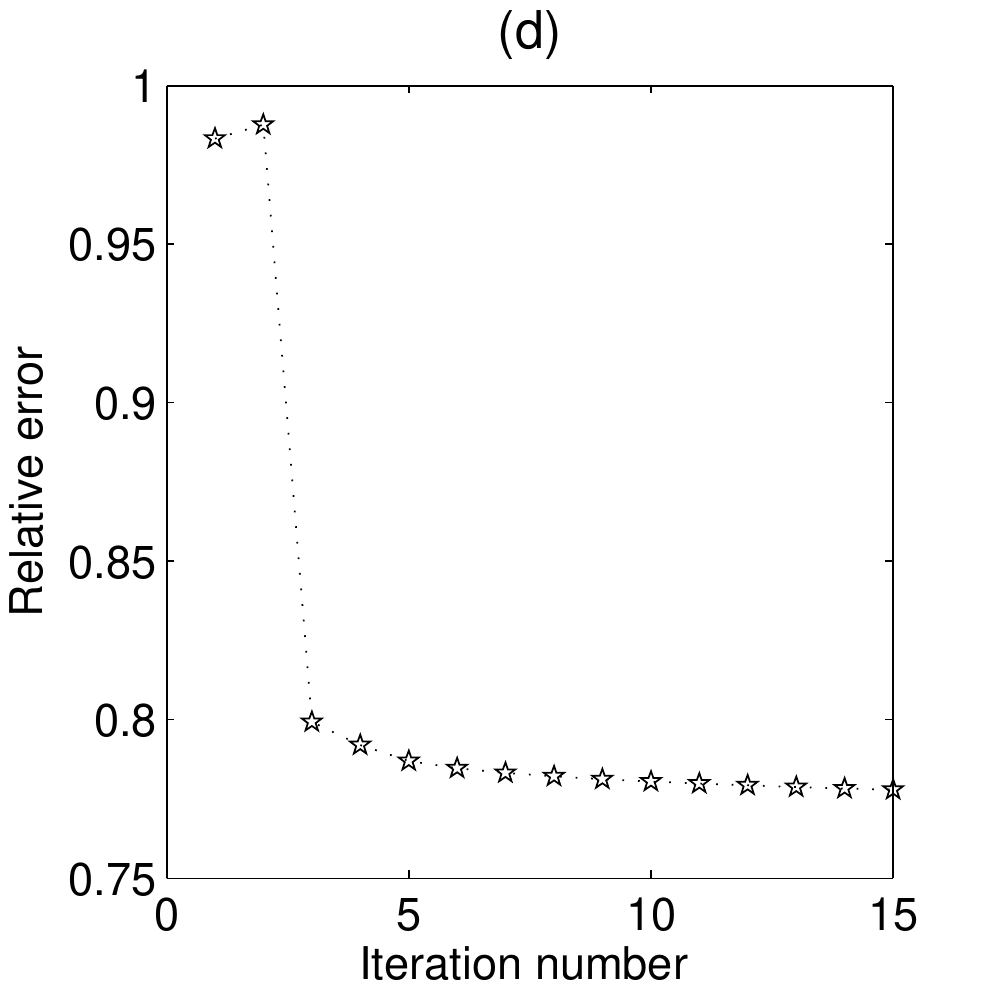}}
\caption {\label{fig5} The results obtained by inverting the data shown in Fig.~\ref{3b} using  the $\chi^2$ principle, the UPRE and the MDP  as the parameter-choice method, respectively. Figs~\ref{4c}, \ref{5c}, \ref{6c}: the progression of the data fidelity $\Phi(\bfd^{(k)})$, the regularization term $\Phi(\bfm^{(k)})$ and the regularization parameter $\alpha^{(k)}$ with iteration $k$ in each case, respectively and in  Figs~\ref{4d}, \ref{5d}, \ref{6d}: the progression of the relative error at each iteration for the same cases. }
\end{figure}
\begin{figure}
\subfigure{\label{5e}\includegraphics[width=.35\textwidth]{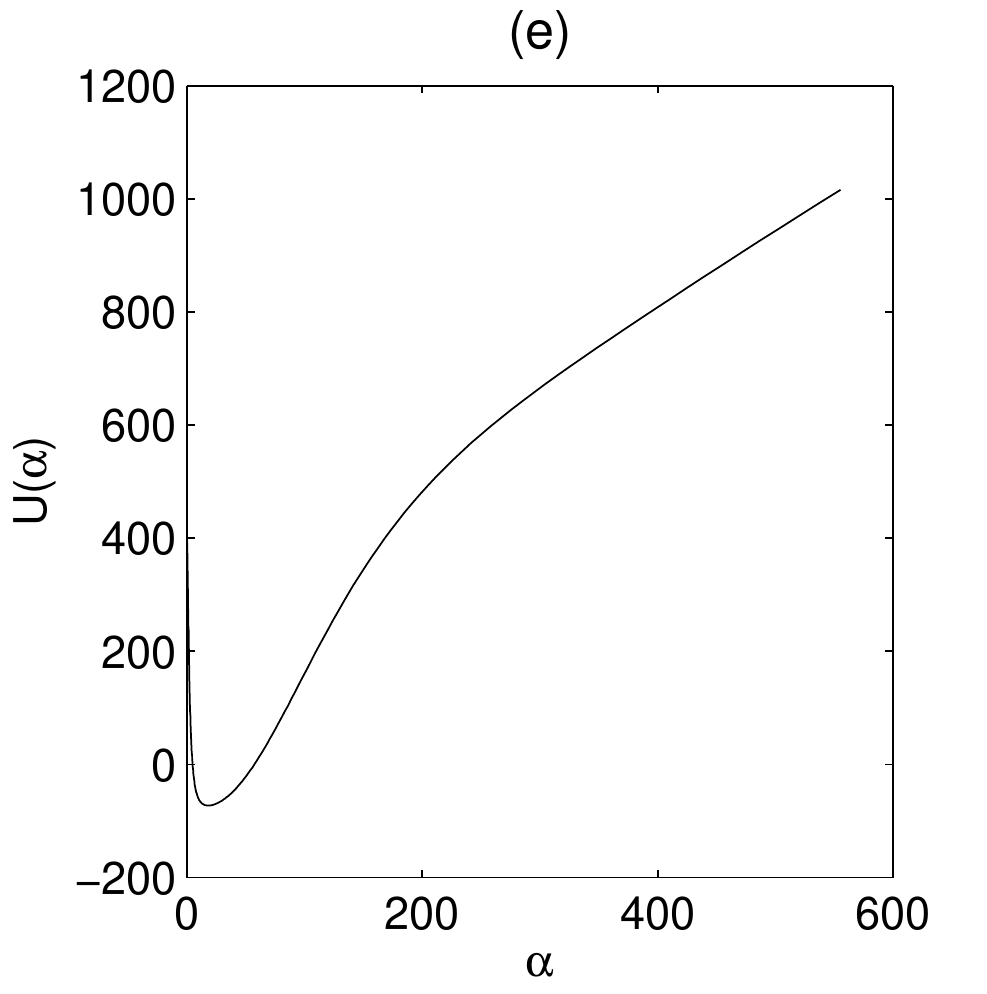}}
\subfigure{\label{5f}\includegraphics[width=.35\textwidth]{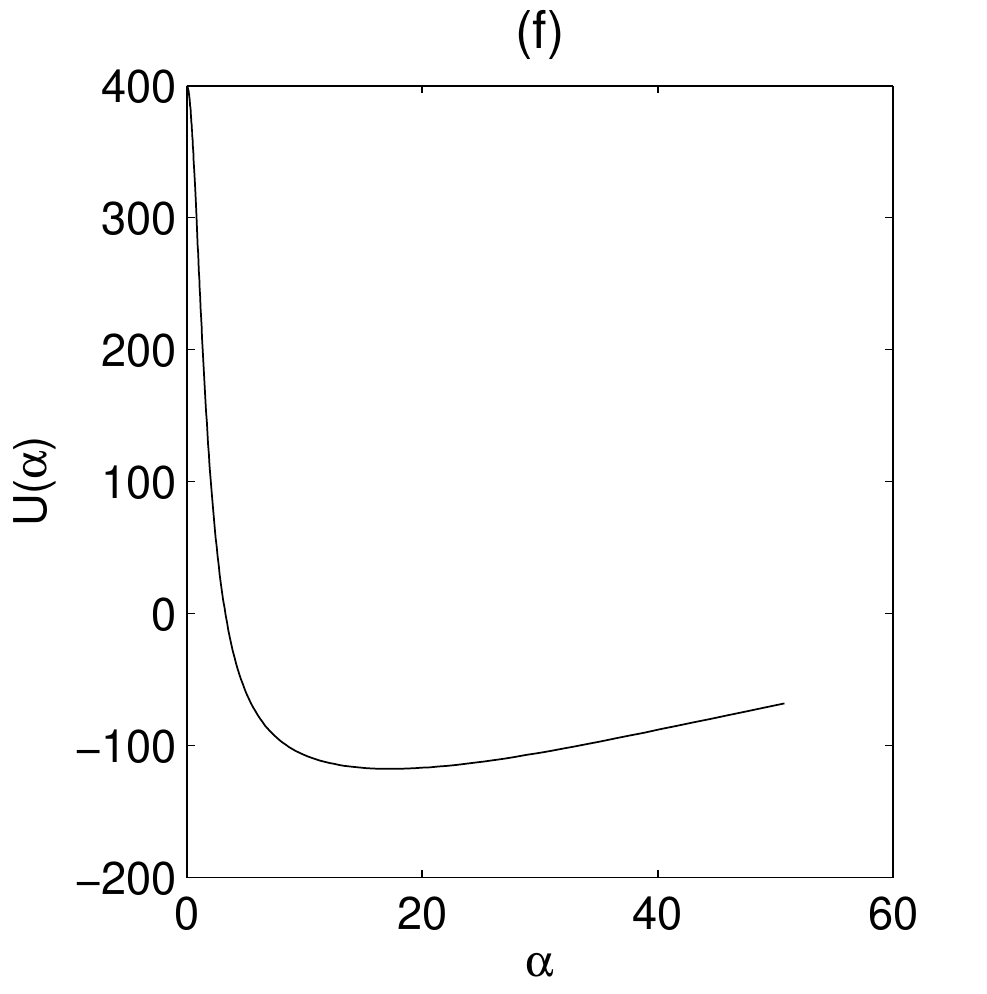}}
\caption {\label{fig6} ; Fig.~\ref{5e}: the UPRE functional at iteration $3$; Fig.~\ref{5f} the UPRE functional at iteration $7$.}
\end{figure}

To illustrate the results summarized in Tables~\ref{tab1}-\ref{tab3},   Figs~\ref{fig4}-\ref{fig6} provide details for a representative case, sample $c=4$ for the second noise level, $(\eta_1, \eta_2)=(0.02, 0.005)$.  Here Figs~\ref{4a}, \ref{5a}, \ref{6a} show the inverted data in cross section at $y=525$~m and Figs~\ref{4b}, \ref{5b}, \ref{6b} the plane sections at $z=100$~m and $z=350$~m. The progression of the data fidelity $\Phi(\bfd^{(k)})$, the regularization term $\Phi(\bfm^{(k)})$ and regularization parameter $\alpha^{(k)}$ with iteration $k$ are presented in Figs~\ref{4c}, \ref{5c}, \ref{6c}, and in  Figs~\ref{4d}, \ref{5d}, \ref{6d} the progression of the relative error. To show that the UPRE functional has a nicely defined minimum we show the functional $U(\alpha)$ at the third and seventh iterations in Figs~\ref{5e}-\ref{5f}. 
 In all cases the algorithms produce a dramatic decrease in the relative error by the third iteration, after which the error decreases monotonically, but with a slower rate for the MDP.  At the same time the regularization parameter appears to stabilize in each case after the fifth iteration, which is contrary to what one would see by  using  iterated Tikhonov, which forces the parameter slowly to zero, e.g. \cite{Tikhonov,Zhdanov}.  The stabilization observed here suggests that it may be sufficient to carry out the regularization parameter estimation only for a limited number of initial steps, but would  require introduction of yet another parameter to assess for stabilization of $\alpha$. Moreover, further experiments not reported here demonstrate that a dramatic increase in iterations is possible for $\alpha^{(k)}$ not chosen to represent the error levels in the current iteration. Thus, it is important to continue to update $\alpha$ every step of the iteration.

\subsection{Synthetic example: Cube}\label{cube}
\begin{figure}
\subfigure{\label{7a}\includegraphics[width=.35\textwidth]{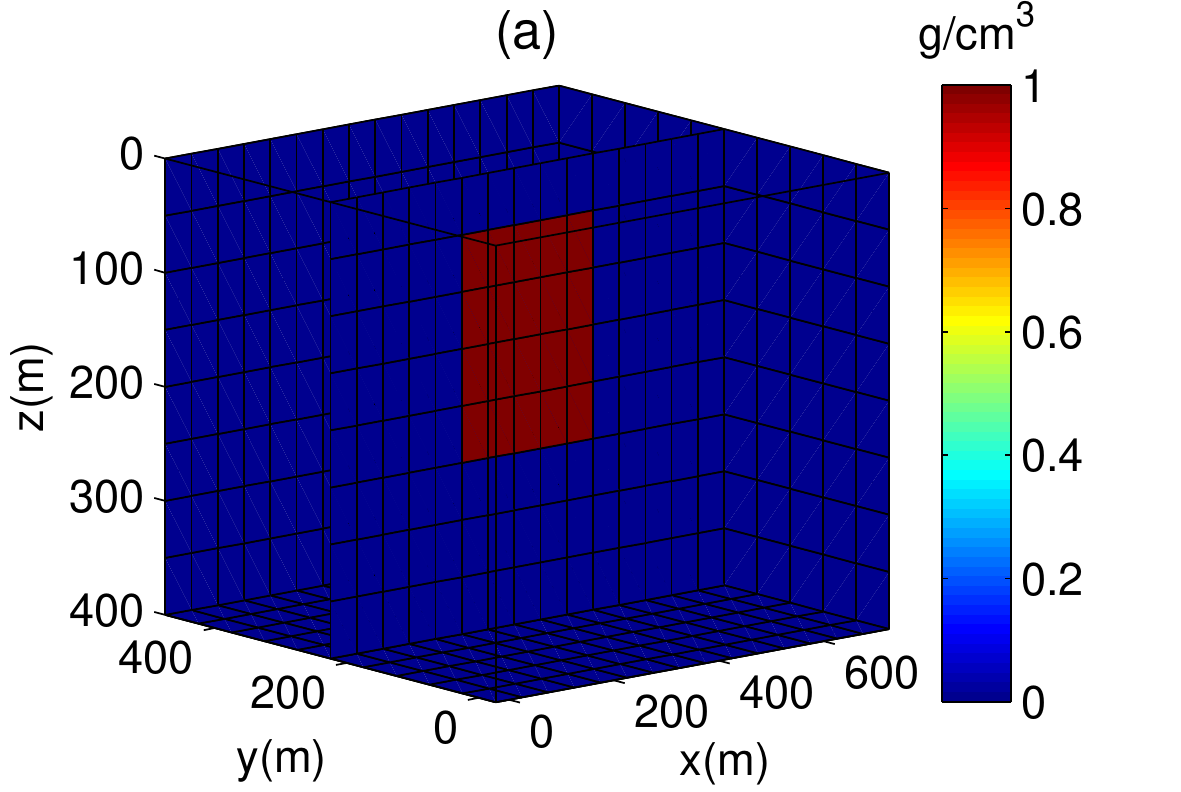}}
\subfigure{\label{7b}\includegraphics[width=.35\textwidth]{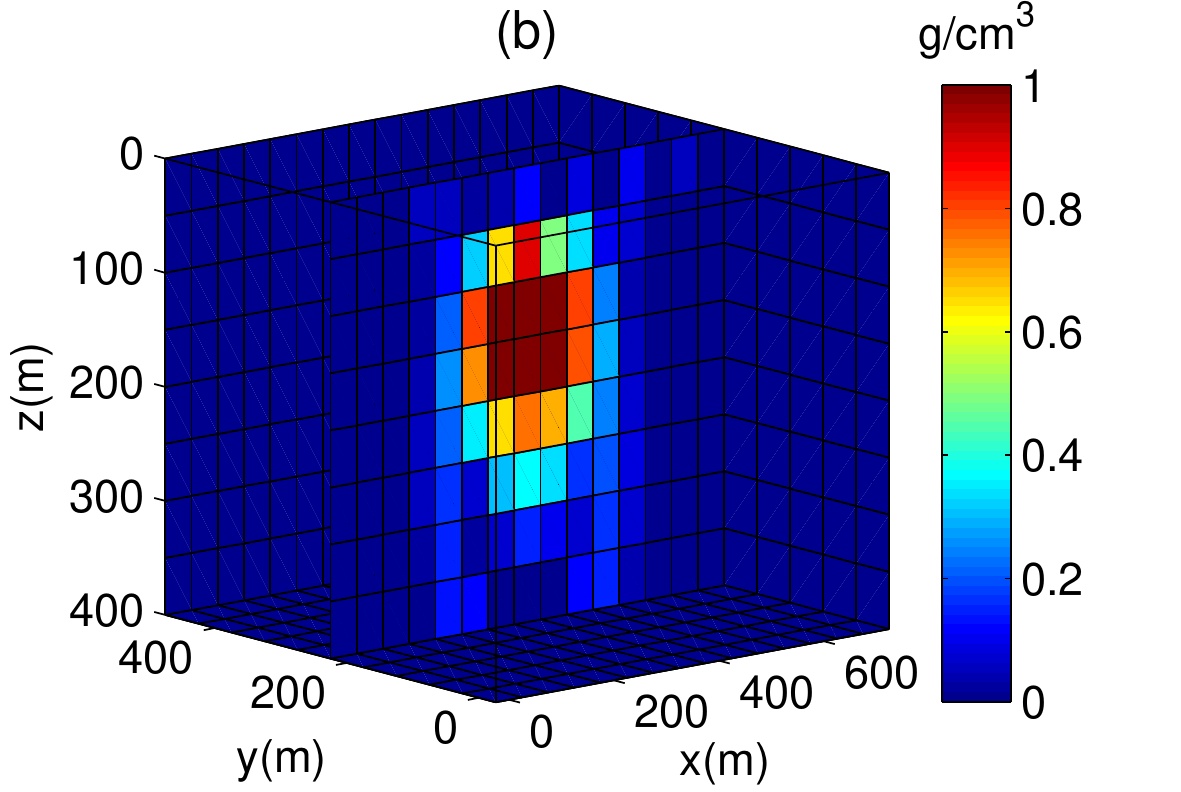}}\\
\subfigure{\label{7c}\includegraphics[width=.35\textwidth]{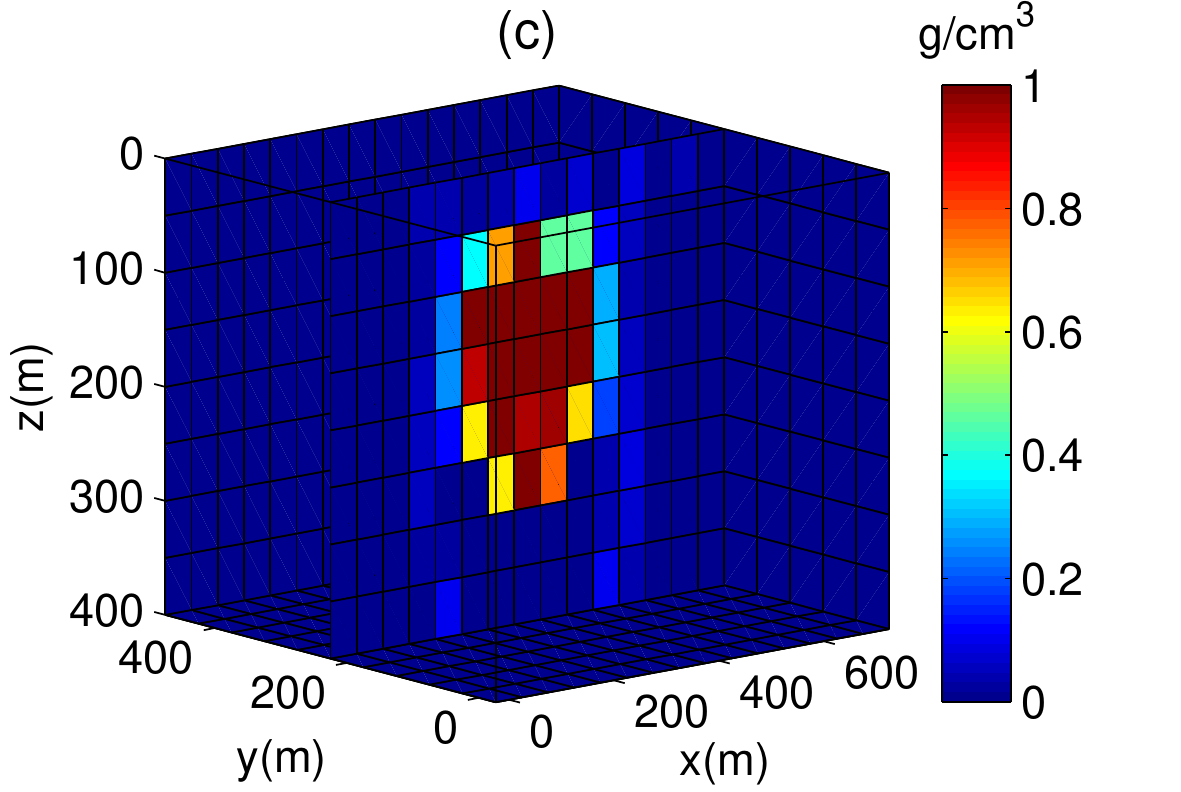}}
\subfigure{\label{7d}\includegraphics[width=.35\textwidth]{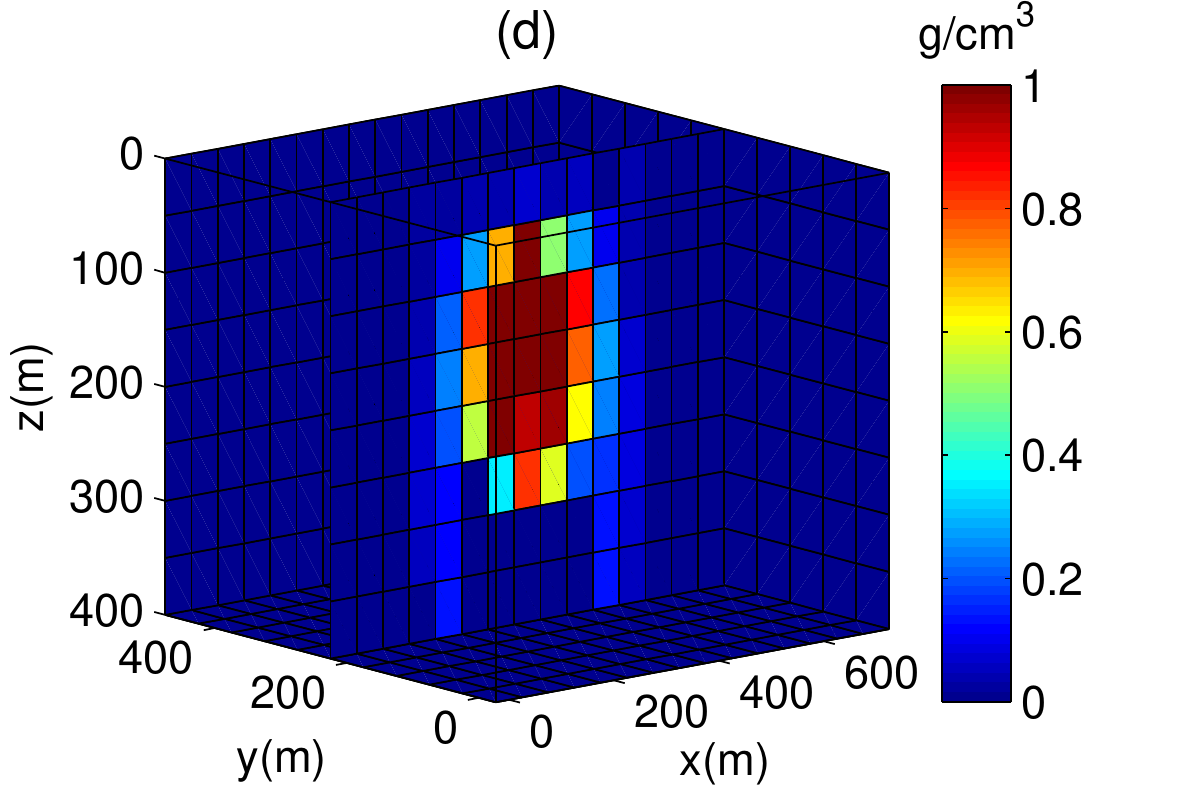}}
\caption{\label{fig7} Fig~\ref{7a}:  Model of a cube  on  an homogeneous background. The density contrast of the cube is $1$~g$/$cm$^3$. Fig.~\ref{7b}: The density model obtained using the $\chi^2$ principle; Fig.~\ref{7c}: The density model obtained using the UPRE; Fig.~\ref{7d}:  The density model obtained using the MDP.}
\end{figure}

As a second example we choose a cube with dimension $250$~m$\times~200$~m~$\times~200$~m with  density  contrast $1$~g$/$cm$^3$ on  an homogeneous background, Fig.~\ref{7a}. Simulation data, $\bfd$, are calculated over a $15$ by $10$ grid with spacing $\Delta=50$~m  on the surface, using the same  three noise levels  as for the dike simulations. For inversion the subsurface is divided into $ 15 \times 10 \times 8 = 1200$ cells  each of size $\Delta= 50$~m.  The simulations are set up as for the case of the dike  and the results of the inversions are summarized in  Tables ~\ref{tab4} - \ref{tab6}, for parameter estimation using the $\chi^2$ principle, the UPRE method, and the MDP method, respectively.  An illustration of these results is given in Fig.~\ref{fig7} for the case $c=5$ for noise level three, $(\eta_1,  \eta_2)=(0.03, 0.01)$. These results corroborate the conclusions about the performance of each method for the dike simulations.

\begin{table}
\caption{The inversion results obtained by inverting the data from the cube contaminated with the  first noise level, $(\eta_1,  \eta_2)=(0.01, 0.001))$, average(standard deviation)  over $10$ runs.}\label{tab4}
\begin{tabular}{c  c  c  c  c }
\hline
Method&     $\alpha^{(1)},\gamma=1.5$& $\alpha^{(K)}$& Relative error& Number of iterations \\ \hline
$\chi^2$ principle&  1662& 98.4(21.8)& 0.4144(0.0058)& 4.9(0.8)\\
UPRE&  1662& 43.6(3.9)& 0.4150(0.0055)& 4.3(0.5)\\
MDP&  1662& 107(4.3)& 0.4225(0.0050)& 8.1(0.33)\\ \hline
\end{tabular}
\end{table}

\begin{table}
\caption{The inversion results obtained by inverting the data from the cube contaminated with the  second noise level, $(\eta_1, \eta_2)=(0.02, 0.005)$, average(standard deviation)  over $10$ runs.}\label{tab5}
\begin{tabular}{c  c  c  c  c }
\hline
Method&     $\alpha^{(1)},\gamma=1.5$& $\alpha^{(K)}$& Relative error& Number of iterations \\ \hline
$\chi^2$ principle&  1688& 37.7(5.0)& 0.4200(0.0105)& 5.3(1.3)\\
UPRE&  1688& 18.2(3.0)& 0.4225(0.0196)& 4.9(0.9)\\
MDP& 1688& 36.9(3.6)& 0.4202(0.0198)& 12.0(2.3)\\ \hline
\end{tabular}
\end{table}

\begin{table}
\caption{The inversion results obtained by inverting the data from the cube contaminated with the  third noise level, $(\eta_1,  \eta_2)=(0.03, 0.01)$, average(standard deviation)  over $10$ runs.}\label{tab6}
\begin{tabular}{c  c  c  c  c }
\hline
Method&     $\alpha^{(1)},\gamma=1.5$& $\alpha^{(K)}$& Relative error& Number of iterations \\ \hline
$\chi^2$ principle&  1699& 65.4(24.8)& 0.4878(0.0324)& 4.1(0.33)\\
UPRE&  1699& 16.7(2.8)& 0.4769(0.0397)& 4.1(0.6)\\
MDP&  1699& 23.8(6.6)& 0.4808(0.0305)& 5.9(1.2)\\ \hline
\end{tabular}
\end{table}

\subsection{Solution by the generalized singular value decomposition}
In prior work we have used the GSVD to find  $\bfz(\alpha)$ in \eqref{tikhonov3} in place of the SVD as used for the results presented in Sections~\ref{dike}-\ref{cube}. Here we are not presenting the results using the GSVD. There is no difference in the conclusions that may be deduced concerning the efficacy of the regularization parameter estimators but the GSVD is noticeably more expensive. Indeed there is no difference in the results, i.e. $\alpha^{(K)}$, $K$ and the relative errors are the same, but for a greater computational cost, in our implementation the GSVD algorithm is about  $30\%$ more expensive to run. 
In particular, we note that the standard algorithms for finding a GSVD, first find the SVD of the system matrix $\tilde{G}$. On the other hand, for the implementation using the SVD for $\tilde{\tilde{G}}$ one needs only the SVD and the calculation of the inverse for matrix $D$ which in this case is trivially obtained noting that $D$ is diagonal. It is thus not surprising to find that it is more efficient to use the SVD in place of the GSVD.

\section{Real data}\label{real}
\subsection{Geological context}\label{geology}
The field data which is used for modeling are acquired over an area located in  the south-west of Iran where a dam, called Gotvand, is constructed on the Karoon river. Tertiary deposits of the Gachsaran formation are the dominant geological structure in the area. It is mainly comprised of  marl, gypsum, anhydrite and halite. There are several solution cavities in the halite member of the Gachsaran formation which have outcropped with sink-holes in the area. One of the biggest sink-holes is located in the south-eastern part of the survey area and is called the Boostani sink-hole. The main concern is that it is possible that cavities at the location of the Boostani sink-hole may be connected to several other cavities toward the west and the north and joined to the Karoon river. This can cause a serious leakage of water after construction of the dam or may cause severe damage to the foundations of  the dam.

\subsection{Residual Anomaly}\label{residualA}
\begin{figure}
\includegraphics[width=.8\textwidth]{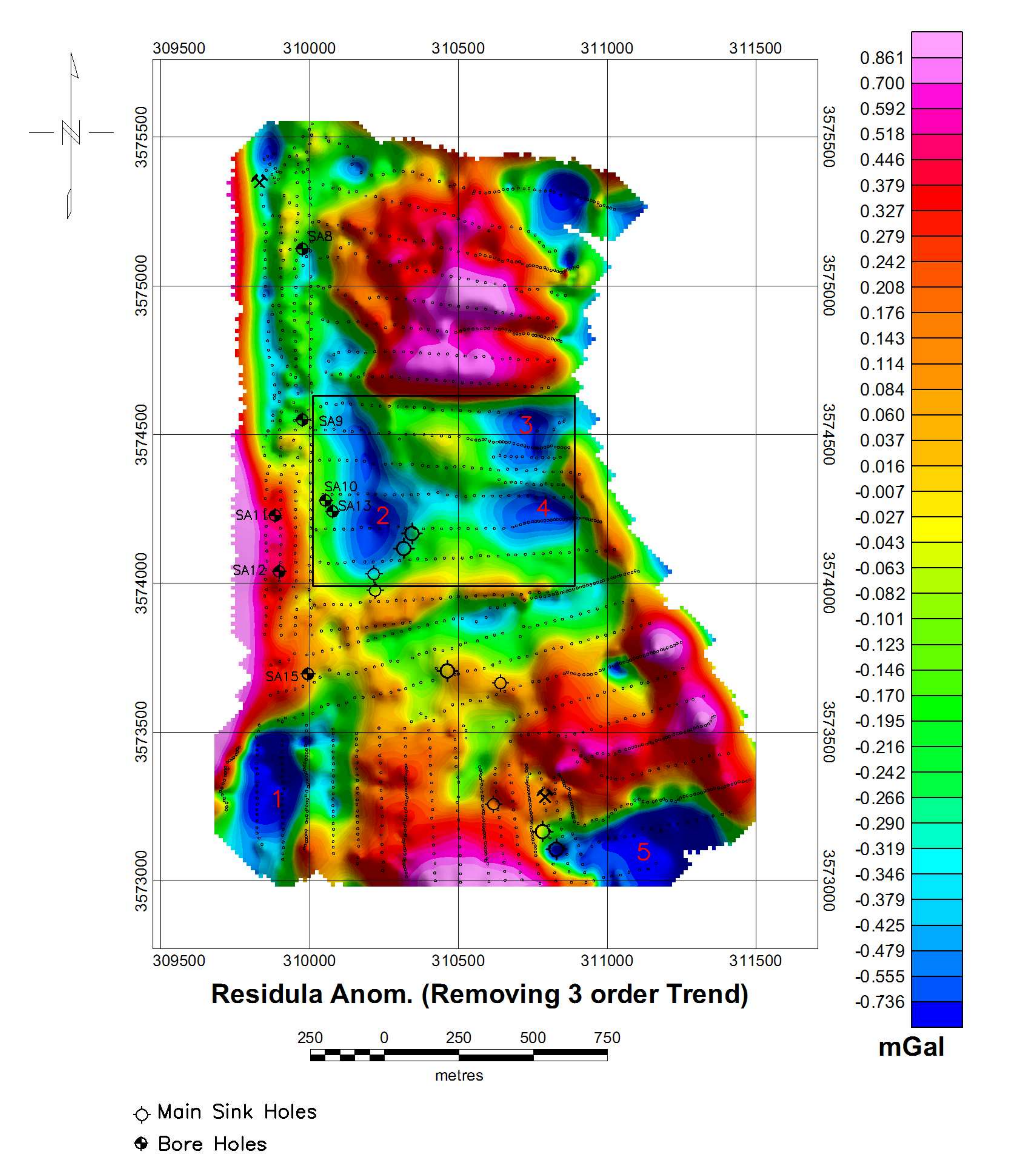}
\caption {\label{fig8} Residual anomaly map over the Gotvand dam site.}
\end{figure}

\begin{figure}
\includegraphics[width=.8\textwidth]{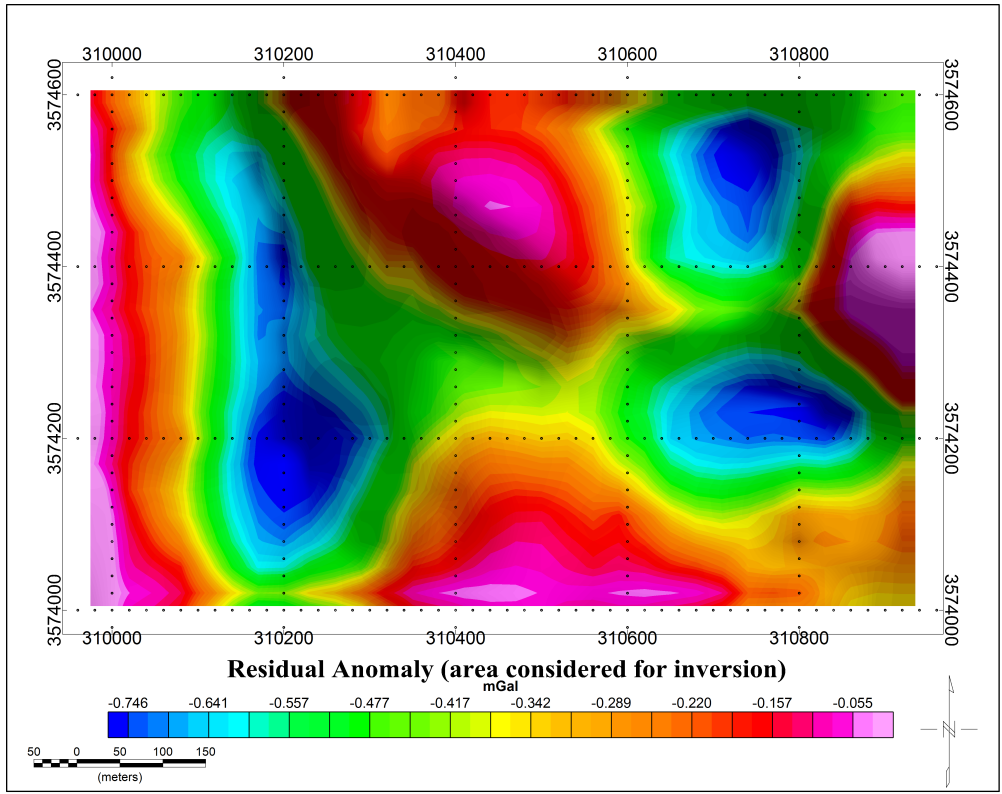}
\caption {\label{fig9} Residual anomaly selected for inversion.}
\end{figure}

The gravity measurements were undertaken by the gravity branch of the Institute of Geophysics, Tehran University. Measurements were taken at  $1600$ stations such that separation between points along the profiles is about $10$~m and separation between profiles is $30$~m to $50$~m. Data were corrected for effects caused by variation in elevation, latitude and topography to yield the Bouguer gravity anomaly. The residual gravity anomaly has been computed using a  polynomial fitting method, Fig.~\ref{fig8}. The six main negative anomalies representing low-density zones are identified on this map. Anomaly $5$ is over the Boostani sink-hole. We have selected  a box including anomalies $2$, $3$ and $4$ for application of the inversion code, Fig.~\ref{fig9}. More details about field procedures, gravity correction and interpretation of the data  are provided in  \cite{Ardestani:2013}.   

\subsection{Inversion results}\label{realresults}
The residual anomaly,  Fig.~\ref{fig9},  was sampled every $30$~m  yielding a box of  $32 \times 20 =640$ gravity points. We suppose that the data is contaminated by error as in the case of the simulations using the noise level case two, $(\eta_1, \eta_2)=(.02, .005)$.  The subsurface is divided into $ 32 \times 20 \times 10 =6400$ cells of size  $\Delta=30$~m in each dimension. Based on geological information a background density
 $2.4$~g$/$cm$^3$ is selected for the inversion and  density is limited by    $\rho_{\mathrm{min}} = 1.5$~g$/$cm$^3$ and $\rho_{\mathrm{max}}=   2.4$~g$/$cm$^3$. The results obtained using all three parameter choice methods are collated in Table~\ref{tab7}.  As for the simulated cases, we find that the final $\alpha$ is larger for both the MDP and $\chi^2$ approaches, suggesting greater smoothing in the solutions. In contrast to the simulated cases, the UPRE requires more iterations to converge, as can be seen in Figs~\ref{10c}, \ref{11c}, \ref{12c}, which show  the progression of the data fidelity $\Phi(\bfd^{(k)})$, the regularization term $\Phi(\bfm^{(k)})$ and the  regularization parameter $\alpha^{(k)}$ with iteration $k$. 
 We stress that the total time for the implementation using the $\chi^2$ principle is about one third of that for the other two methods, requiring in our implementation about $15$ minutes as compared to roughly $40$ minutes. 
 
\begin{figure}
\subfigure{\label{10a}\includegraphics[width=.45\textwidth]{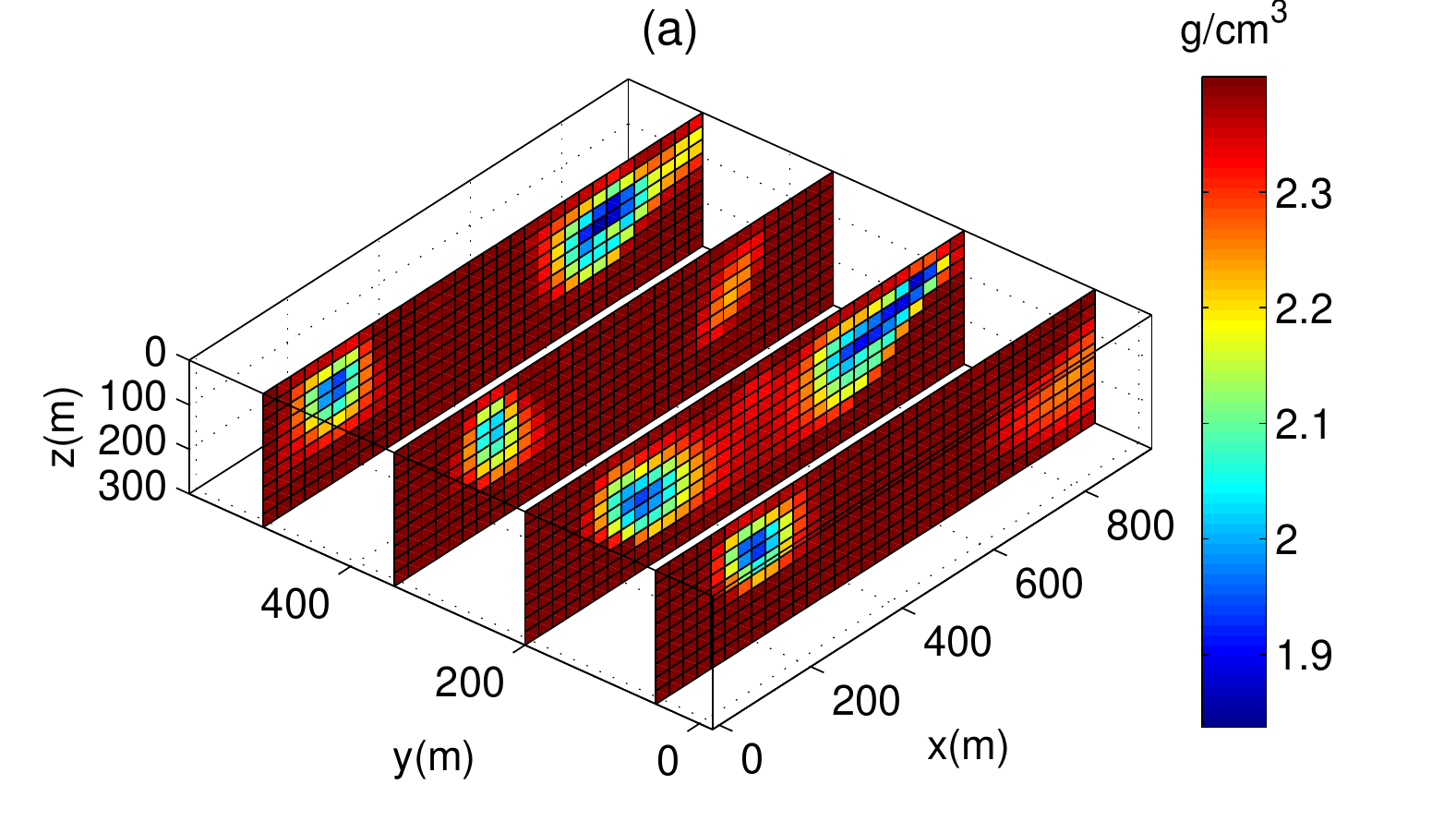}}
\subfigure{\label{10b}\includegraphics[width=.45\textwidth]{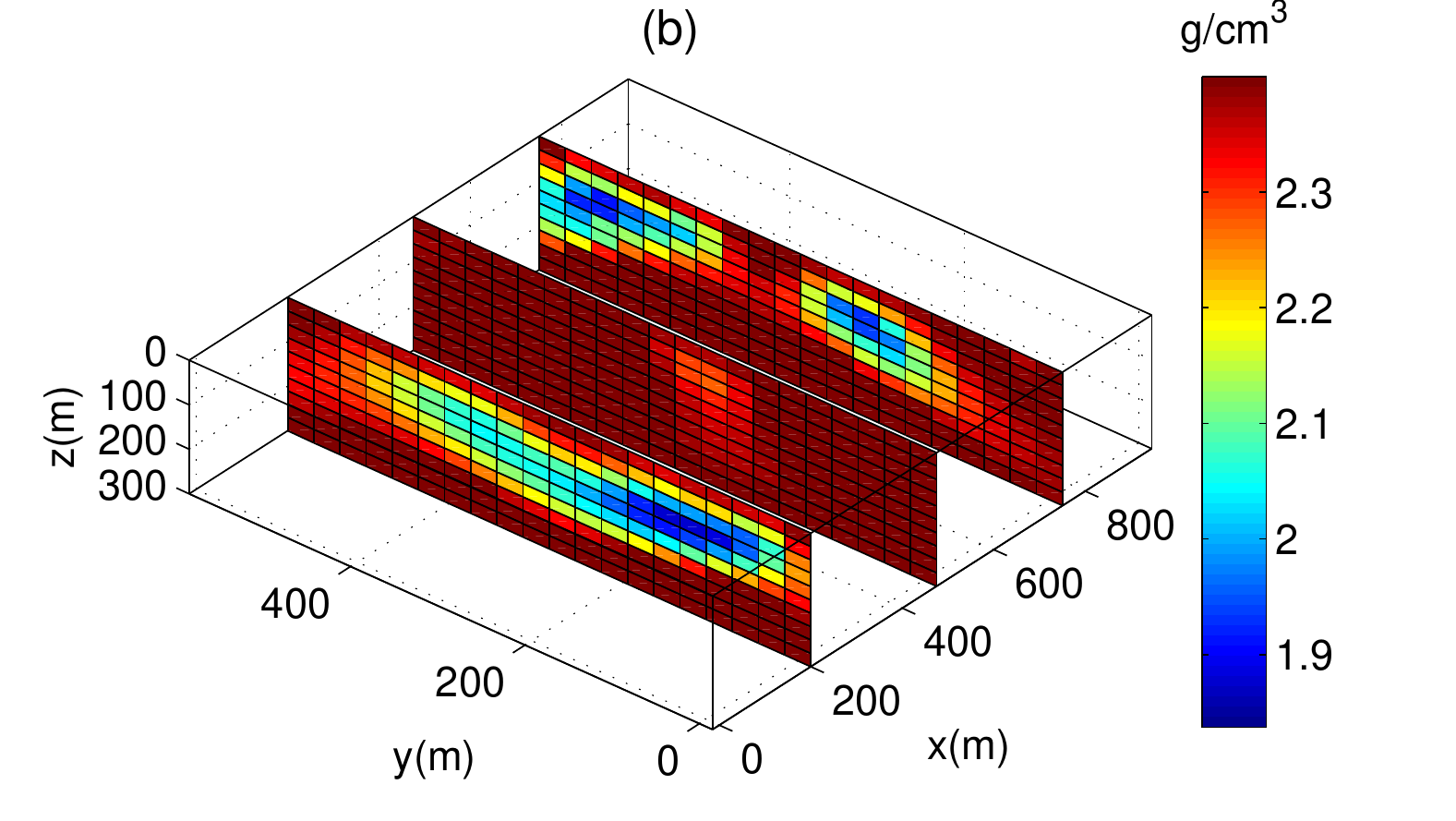}}\\
\subfigure{\label{11a}\includegraphics[width=.45\textwidth]{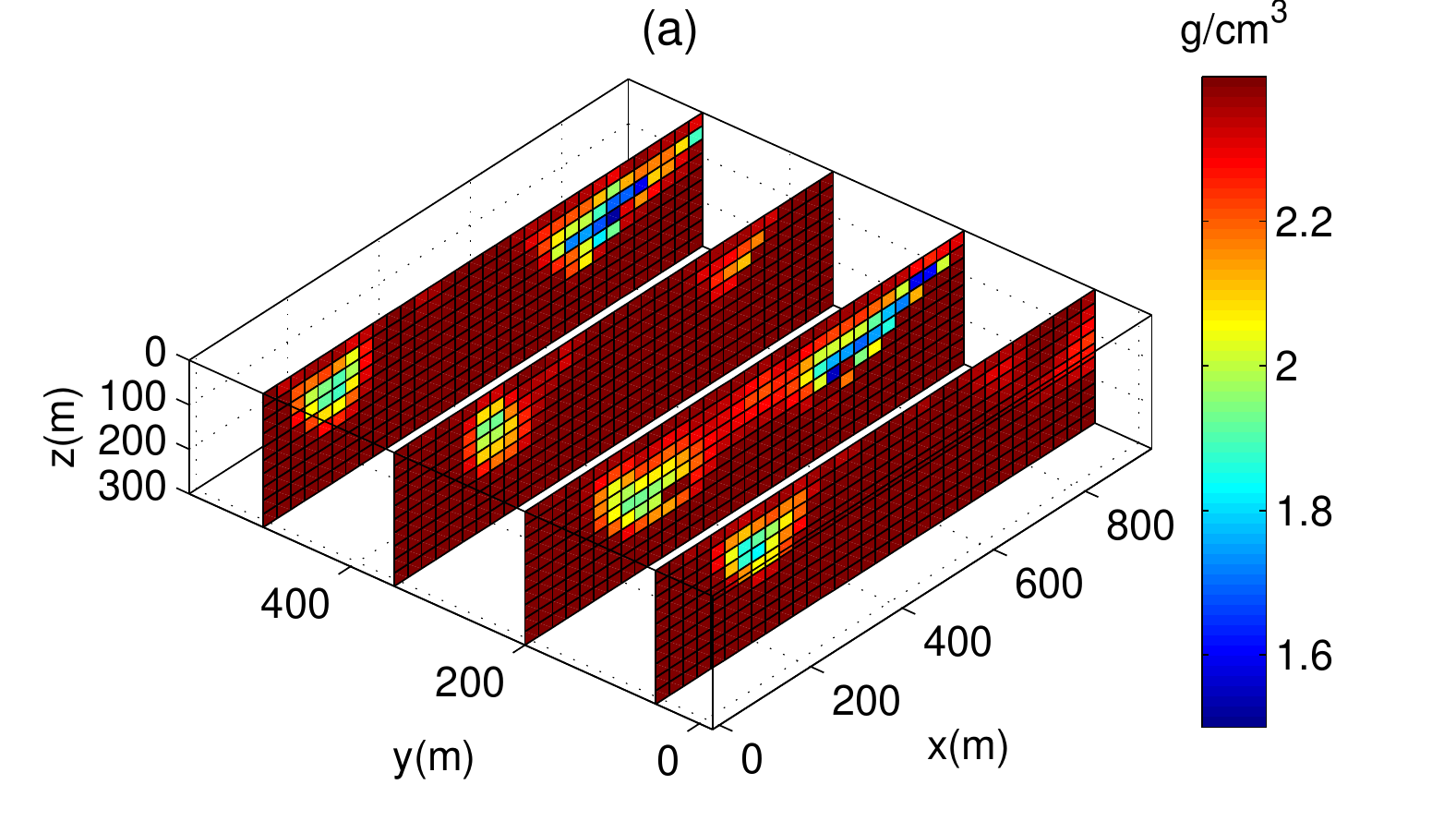}}
\subfigure{\label{11b}\includegraphics[width=.45\textwidth]{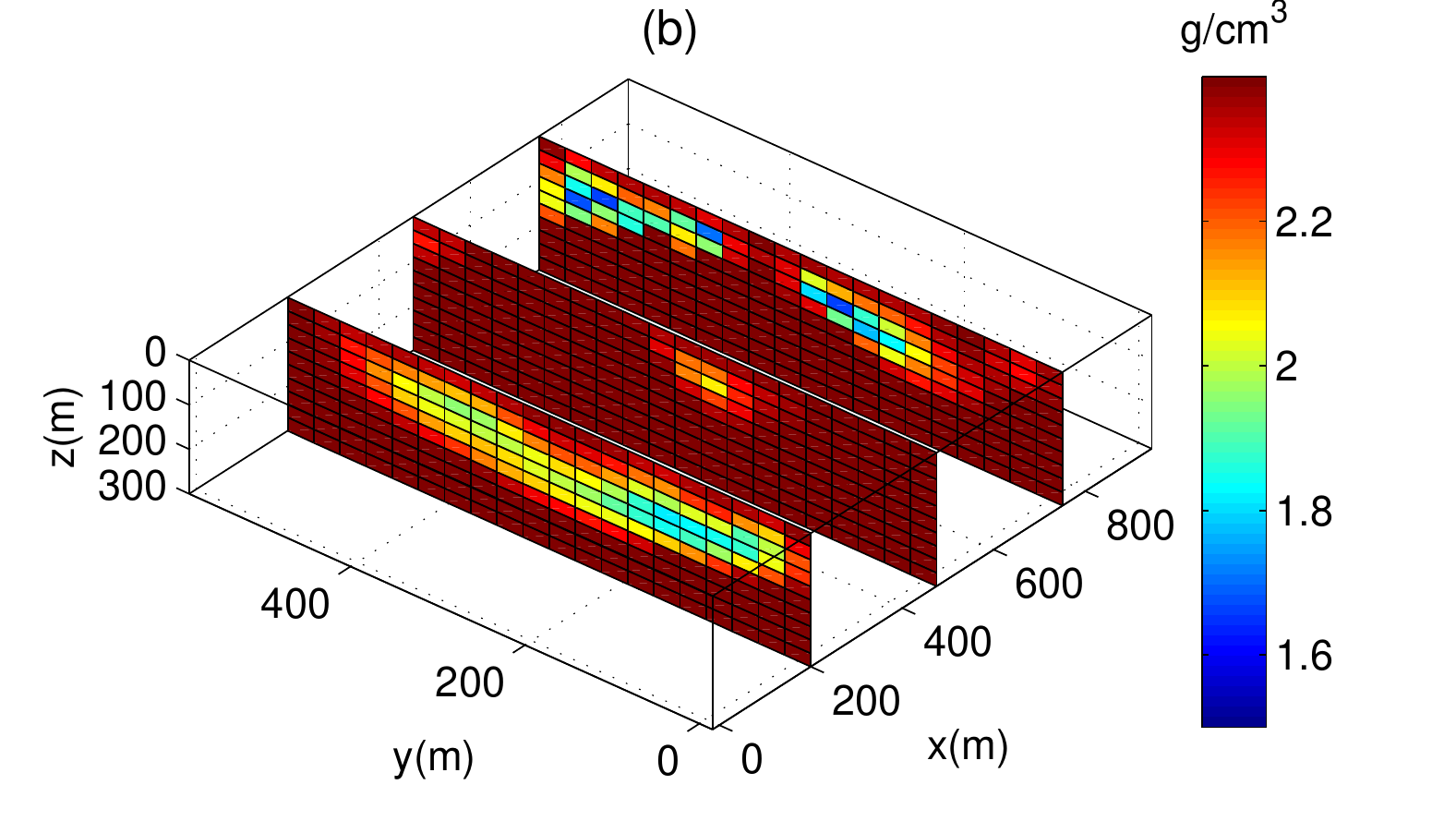}}\\
\subfigure{\label{12a}\includegraphics[width=.45\textwidth]{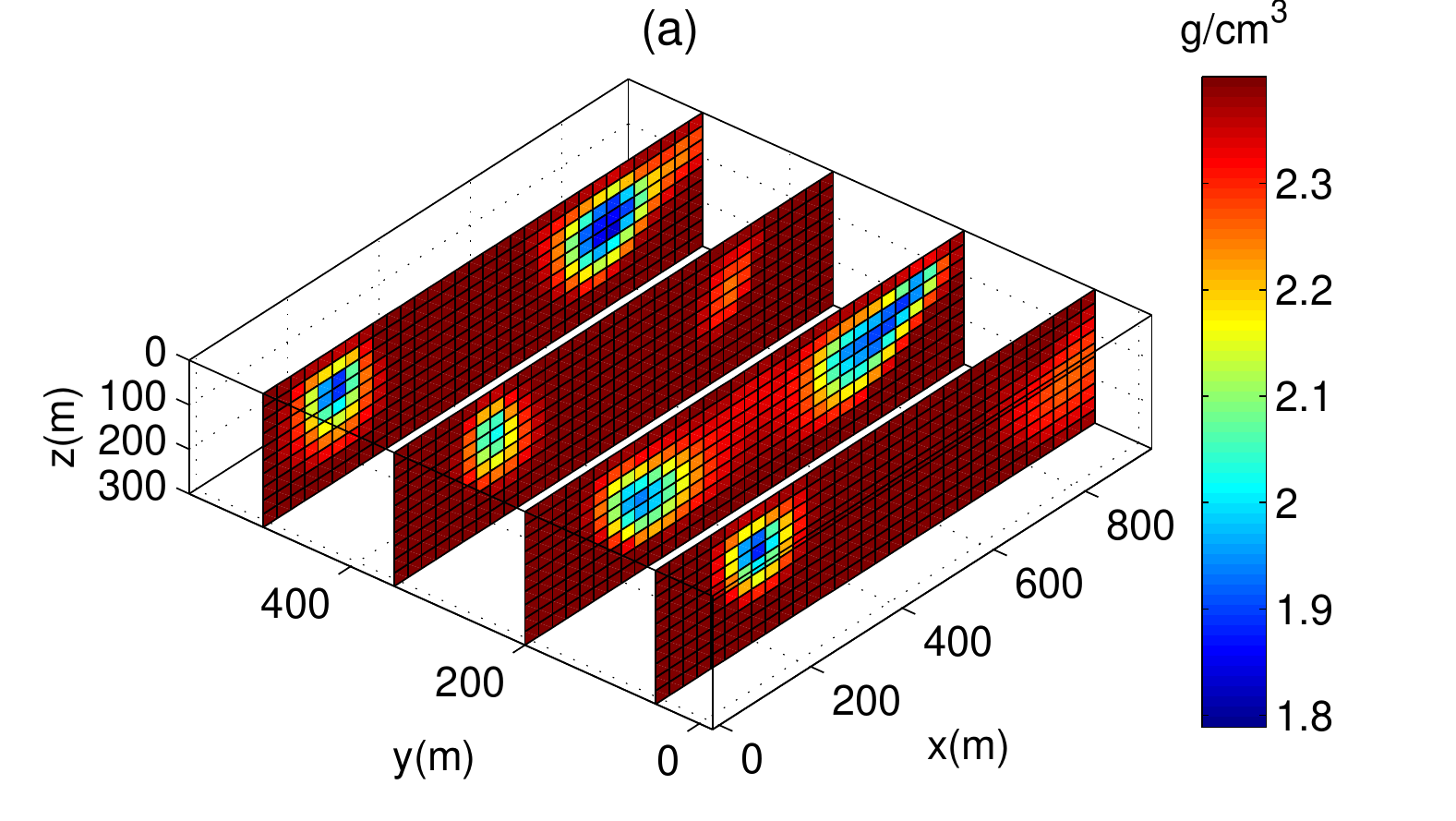}}
\subfigure{\label{12b}\includegraphics[width=.45\textwidth]{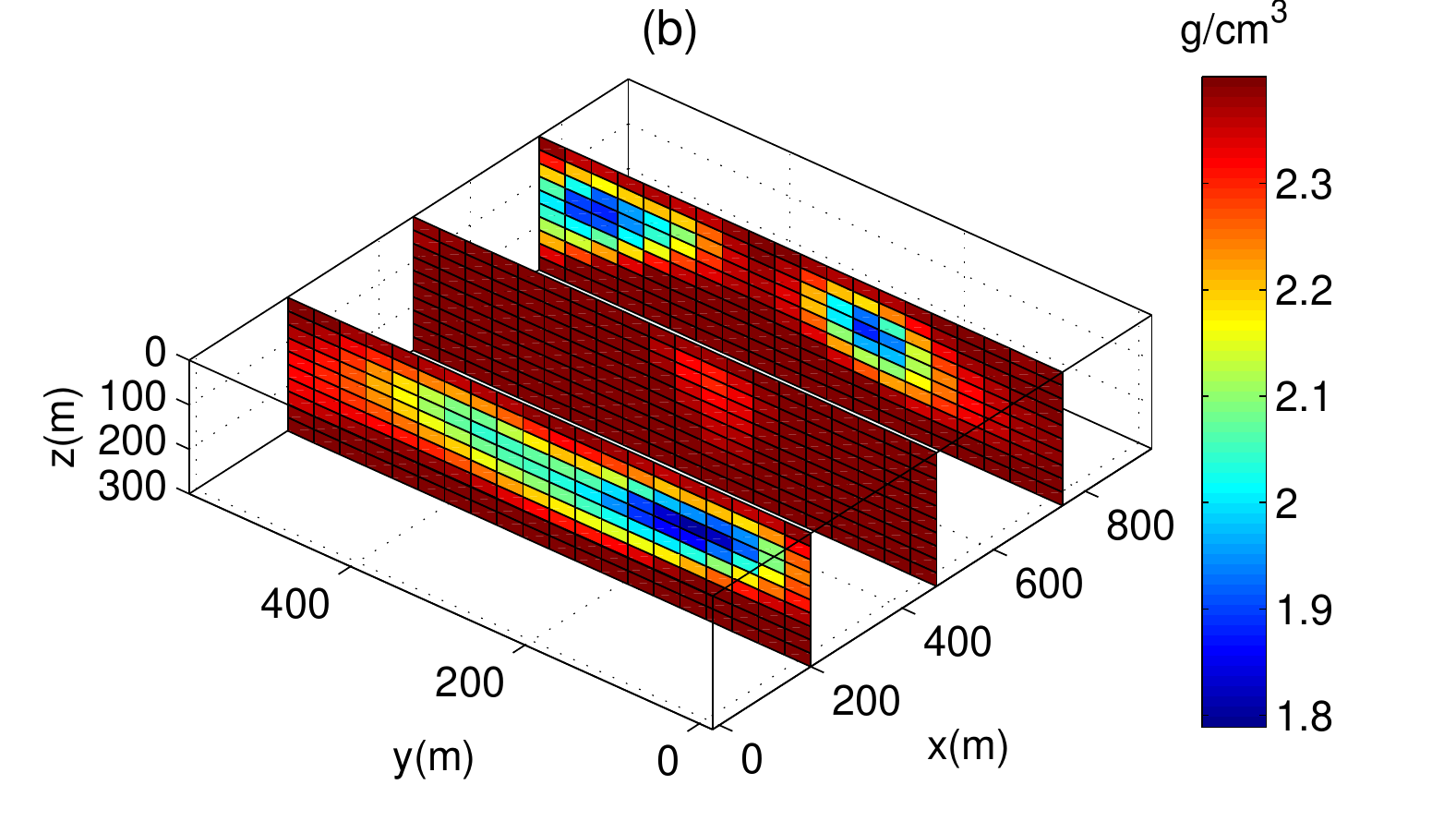}}
\caption {\label{fig10} 
The results obtained by inverting the data shown in Fig.~\ref{fig9} using  the $\chi^2$ principle, the UPRE and the MDP  as the parameter-choice method, respectively. Figs~\ref{10a}, \ref{11a}, \ref{12a}: cross-sections in the  $y-z$ plane in each case, respectively and in  Figs~\ref{10b}, \ref{11b}, \ref{12b}: cross-sections in the $x-z$ plane for the same cases.}
\end{figure}

\begin{figure}
\subfigure{\label{10c}\includegraphics[width=.3\textwidth]{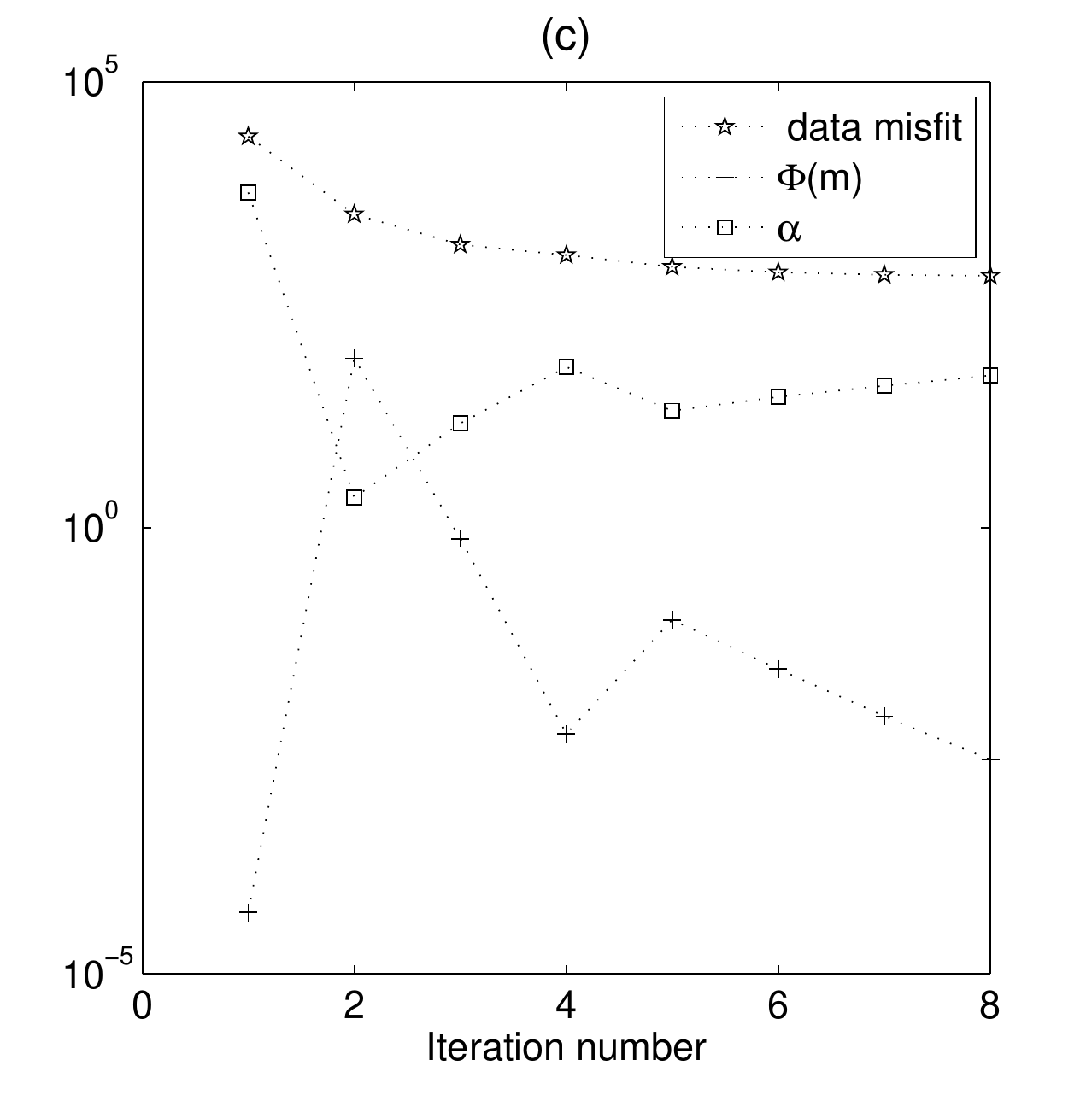}}
\subfigure{\label{11c}\includegraphics[width=.3\textwidth]{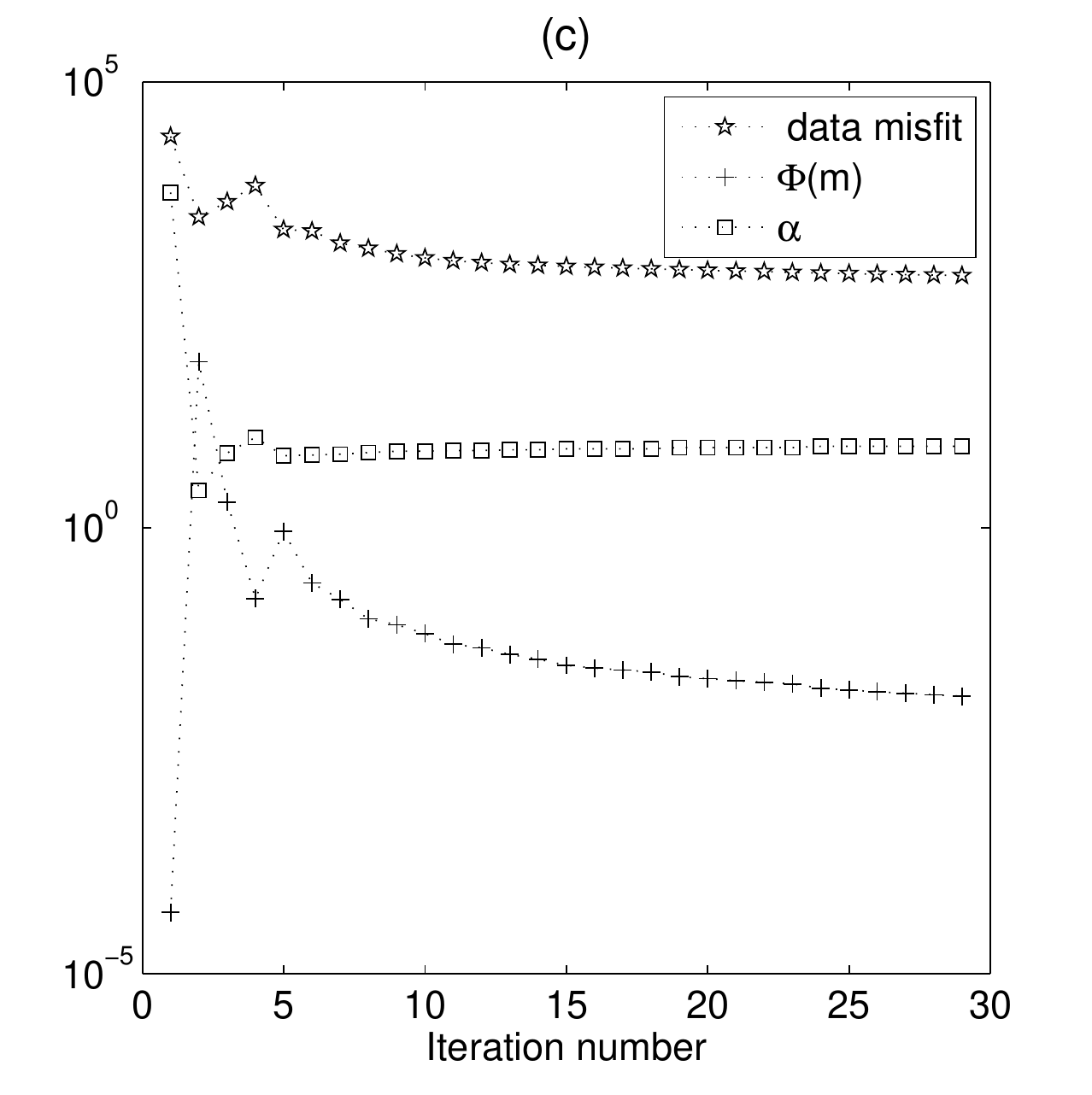}}
\subfigure{\label{12c}\includegraphics[width=.3\textwidth]{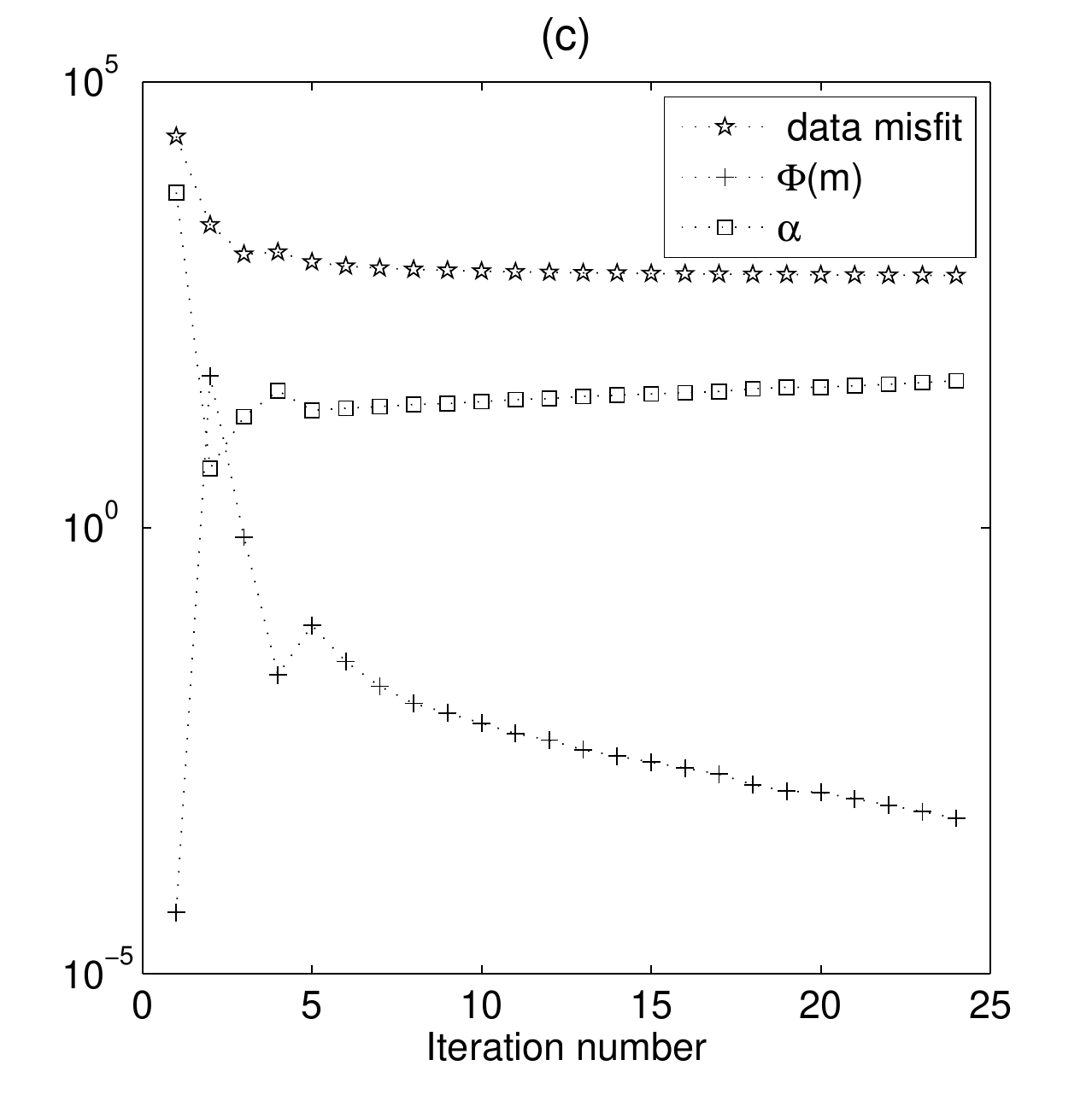}}
\caption{\label{fig11} 
The results obtained by inverting the data shown in Fig.~\ref{fig9} using  the $\chi^2$ principle, the UPRE and the MDP  as the parameter-choice method, respectively. Figs~\ref{10c}, \ref{11c}, \ref{12c}: the progression of the data fidelity $\Phi(\bfd^{(k)})$, the regularization term $\Phi(\bfm^{(k)})$ and the regularization parameter $\alpha^{(k)}$ with iteration $k$ in each case, respectively.}
\end{figure}

 In assessing these results, it is also useful to consider the visualizations of the solutions,  given in Figs~\ref{10a}, \ref{11a}, \ref{12a}, and \ref{10b}, \ref{11b}, \ref{12b}, for the cross sections in the $y-z$ and $x-z$ planes, respectively. Immediate inspection indicates that the solutions using the MDP and $\chi^2$ approach are quite close, while the UPRE differs. Further assessment of the quality of the solutions makes use of our knowledge of the anomalies, the depths of which have been estimated by 3D modeling and are given in Table~\ref{tab8}. Fig.~\ref{fig8} also shows that there are two bore holes in the area near  anomaly two, for which the range of the low-density zone obtained from these bore-holes is also given in  Table~\ref{tab8}. Estimations of the same measures of these anomalies using the reconstructions are also collated in Table~\ref{tab8}. Now it is clear that indeed the reconstructions using the $\chi^2$ and MDP are very close yielding a  range for the density contrast of the  low-density zones $2$ to $4$ of $ 1.8$ to $2.4$. On the other hand, the obtained depths using the UPRE are  closer to those obtained with the  bore-holes,  and  while the density contrast for anomaly $2$ still lies  in the interval $ 1.8$ to $2.4$, for   anomalies  $3$ and $4$ the range is between $1.5$ and $2.4$. We conclude that the UPRE, although needing now more iterations, is potentially more robust than either of the other methods, but that indeed the $\chi^2$ method can be useful for generating solutions more efficiently, with fewer iterations, and might therefore be used  when efficiency is of the highest concern. 
 

\begin{table}
\caption{Results obtained by inverting the data shown in Fig.~\ref{fig9}.}\label{tab7}
\begin{tabular}{c  c  c  c }
\hline
Method&     $\alpha^{(1)},\gamma=1.5$& $\alpha^{(K)}$& Number of iterations   \\ \hline
$\chi^2$ principle&  5743&          51.3 &       8     \\
UPRE              &  5743&          8.2 &        29      \\
MDP               &  5743&          44.5 &         24    \\ \hline
\end{tabular}
\end{table}

\begin{table}
\caption{Depths obtained using 3D modeling.}\label{tab8}
\begin{tabular}{c c c c c c c c c} \hline
Anomaly & \multicolumn{2}{c}{$\chi^2$}& \multicolumn{2}{c}{UPRE}& \multicolumn{2}{c}{MDP} & \multicolumn{2}{c}{Bore-hole} \\ \hline
 & min & max & min & max & min & max & min & max \\ \hline
2 &  30-60 & 150-180 & 60-90 & 150 & 30-60 & 150-180 & 115-150 &   150-160     \\
3 &  30 & 90-180 & 30 & 90-120 & 30 & 90-180 & - &   -     \\
4 &   30 & 150 & 30 & 90 & 30 & 150 &- &   -     \\ \hline
\end{tabular}
\end{table}

\section{Conclusions}\label{conclusion}
The $\chi^2$ and UPRE parameter-choice methods have been introduced in the context of 3D gravity modeling. Presented results validate that both methods are more effective  than the more often used MDP. While the  $\chi^2$ technique is itself very fast for each iteration, requiring only an effective one dimensional root finding algorithm, it also  converges quickly. Thus it is definitely to be preferred over the MDP. On the other hand, the UPRE generally provides results with the least relative error in contrast to the MDP and $\chi^2$ methods, particularly for situations with higher noise levels, even if the results for practical data demonstrate that the number of iterations may be increased. In terms of the implementation of the UPRE, the only disadvantage is that finding the optimal $\alpha$ at each step requires the calculation of the $U(\alpha)$  for a range of $\alpha$. Still we have seen that the minimum of  $U(\alpha)$ is well-defined during the iterations. 

In these results we have presented an algorithm for finding the minimum of the Tikhonov functional using the SVD for the system matrix in standard form \cite{Hansen} at each iteration in contrast to the use of the GSVD for the augmented matrix formed from the system and stabilizing matrices. 
The resulting algorithm is much faster and less memory intense, representing generally $30\%$ savings in our implementation.  Moreover, it has been successfully validated for the modeling of the subsurface for the Gotvand dam site located in south-west  Iran. These results indicate that the low-density zones extend between $60$ and $150$~m in depth, which is in general agreement with measurements obtained from bore-holes. 

While the results here have demonstrated the practicality of the regularization parameter estimation techniques in conjunction with the minimum support stabilizer and the singular value decomposition for 3D focusing gravity inversion, the computational cost per reconstruction is still relatively high. For future work we plan to investigate projected Krylov methods to solve the systems at each iteration. Replacement of the SVD at each step by an iterative technique is straightforward, but the question of determining the optimal regularization parameter for the solution on the underlying Krylov subspace each step is still an unresolved question and worthy of further study for reducing the cost of 3D inversions in complex environments, as well as for inclusion of alternative edge preserving regularizers. 

\section*{Acknowledgments}
Rosemary Renaut acknowledges the support of AFOSR grant 025717: ``Development and Analysis of Non-Classical Numerical Approximation Methods", and 
NSF grant  DMS 1216559:   ``Novel Numerical Approximation Techniques for Non-Standard Sampling Regimes". 

\appendix
\section{The singular value decomposition}\label{appSVD}

The solution of the  regularized problem defined by right preconditioned matrix $\tilde{\tilde G}$ uses the singular value decomposition (SVD) of the matrix $ \tilde{\tilde G}$ . Matrix $\tilde{\tilde G} \in \mathcal{R}^{m \times n}$, $ m < n$, is factorized as $ \tilde{\tilde G} = U \Sigma V^T$.  The singular values are ordered $\sigma_1 \geq \sigma_2 \geq \cdots \geq \sigma_m > 0$ and occur on the diagonal of  $\Sigma \in \mathcal{R}^{m \times n}$ which has  $n-m$ zero columns,  \cite{GoLo:96}. Matrices $U \in \mathcal{R}^{ m\times m}$ and $V \in \mathcal{R}^{ n\times n}$ are row and column orthonormal. Then the solution of the regularized problem with parameter $\alpha$ is
\begin{eqnarray}\label{svdsoln}
\bfz(\alpha) &=& \sum_{i=1}^{m} \frac{\sigma^2_i}{\sigma^2_i+\alpha^2} \frac{\bfu^T_{i}\tilde{\bfr}}{\sigma_i} \bfv_{i} = \sum_{i=1}^{m} f_i(\alpha) \frac{s_i}{\sigma_i} \bfv_{i} \quad s_i=\bfu^T_{i}\tilde{\bfr}\\
f_i(\alpha)&=&\frac{\sigma^2_i}{\sigma^2_i+\alpha^2}, \,1\le i \le m,  \quad s_{i}=\bfu_{i}^T \tilde{\bfr},   \label{filterfacs}
\end{eqnarray}
where $\bfu_i$  and  $\bfv_i$ are the $i$th columns of matrices $U$ and $V$ and $f_i(\alpha)$ are the filter factors.

\section{Regularization parameter estimation}\label{regparam}
\subsection{Morozov discrepancy principle}
Using the SVD for $\tilde{\tilde G}$,  the MDP for finding $\alpha$ solves 
\begin{eqnarray}\label{MDPSVD}
\sum_{i=1}^{m} \left( \frac{1}{\sigma_{i}^{2}\alpha^{-2}+1}\right)^2 (\bfu_{i}^T \tilde{\bfr})^2-m = 0.
\end{eqnarray}
\subsection{Unbiased predictive risk estimator}
Regularization parameter $\alpha$ is found to minimize the functional 
\begin{eqnarray}\label{upresvd}
U(\alpha)=\sum_{i=1}^{m} \left( \frac{1}{\sigma_{i}^{2}\alpha^{-2}+1}\right)^2 (\bfu_{i}^T \tilde{\bfr})^2+2\left(\sum_{i=1}^{m}f_{i}(\alpha)\right)-m. 
\end{eqnarray}
\subsection{The $\chi^2$ principle}
Parameter $\alpha$ is found as the root of 
\begin{eqnarray}\label{chi2GSVD}
\sum_{i=1}^{m} \left( \frac{1}{\sigma_{i}^{2}\alpha^{-2}+1}\right) (\bfu_{i}^T \tilde{\bfr})^2-m = 0.
\end{eqnarray}


\end{document}